\documentclass[10pt
	]{amsart}
\pdfoutput=1

\usepackage{graphicx}
\usepackage{amssymb,amsthm,amsmath,mathrsfs,bm,braket,marginnote}
\usepackage{enumerate}
\usepackage{appendix}
\usepackage{extarrows}
\usepackage{multirow}

\usepackage{pgf}
\usepackage{pgfplots}
\usepackage{tikz}
\usetikzlibrary{arrows,calc}
\usepackage{verbatim}
\usetikzlibrary{decorations.pathreplacing,decorations.pathmorphing}
\usepackage[numbers,sort&compress]{natbib}
\usepackage{dsfont}

\usepackage{array,dcolumn}
\usepackage{graphicx}
\usepackage[hyperfootnotes=true,
		pdffitwindow=true,
		plainpages=false,
		pdfpagelabels=true,
		pdfpagemode=UseOutlines,
		pdfpagelayout=SinglePage,
		hyperindex,]{hyperref}
\usepackage[activate={true,nocompatibility},spacing,kerning]{microtype}
\usepackage[utf8]{inputenc}

\newcommand{\mathsym}[1]{{}}
\newcommand{\unicode}[1]{{}}

\theoremstyle{plain}
\newtheorem{theorem}{Theorem}

\newtheorem{corollary}[theorem]{Corollary}
\newtheorem{proposition}[theorem]{Proposition}

\theoremstyle{definition}
\newtheorem{definition}{Definition}
\newtheorem{example}{Example}

\theoremstyle{remark}
\newtheorem{remark}[theorem]{Remark}

\newcommand{\R}{\mathbb R}

\newcommand{\tr}{\mathrm{Tr}}
\newcommand{\diag}{\mathrm{diag}}
\newcommand{\dv}{\mathrm{d}}

\newcommand{\U}{\mathrm U}

\newcommand{\s}{\mathcal S}

\newcommand{\phq}{{_p}H_q}

\numberwithin{equation}{subsection}

\begin{document}
	
\title{Cyclic P\'olya Ensembles on the Unitary Matrices and their Spectral Statistics}
	
\author{Mario Kieburg}
\email[MK]{m.kieburg@unimelb.edu.au}
	
\author{Shi-Hao Li}
\email[SL]{lishihao@lsec.cc.ac.cn, shihao.li@scu.edu.cn}
	
\author{Jiyuan Zhang}
\email[JZ]{jiyuanzhang.ms@gmail.com}
	
\author{Peter J. Forrester}
\email[PJF]{pjforr@unimelb.edu.au}
	
\address[MK, JZ and PJF]{School of Mathematics and Statistics, University of Melbourne, 813 Swanston Street, Parkville, Melbourne VIC 3010, Australia}

\address[SL]{Department of Mathematics, Sichuan University, Chengdu, 610064, China}
\date{\today}

\begin{abstract}
\noindent
A framework to study the eigenvalue probability density function for products of unitary random matrices with an invariance property is developed.
This involves isolating a class of invariant unitary matrices, to be referred to as cyclic P\'olya ensembles, and examining their properties with respect to the 
spherical transform on $\U(N)$.
 Included in the cyclic P\'olya ensemble class are Haar invariant unitary matrices,
the circular Jacobi ensemble, known in relation to the Fisher-Hartwig singularity in the theory of Toeplitz determinants, as well as the heat kernel for Brownian motion on the unitary group. We define cyclic P\'olya frequency functions and show their relation to the cyclic P\'olya ensembles, and give a uniqueness statement for the corresponding weights. The natural appearance of bilateral hypergeometric series is highlighted, with this special function playing the role of the Meijer G-function in the transform theory of unitary invariant product of positive definite matrices. 
We construct a family of functions forming bi-orthonormal pairs which underly the correlation kernel of the corresponding determinantal point processes,
and furthermore obtain an integral formula for the correlation kernel involving just two of these functions.
\end{abstract}

\maketitle

\section{Introduction}\label{sec:intro}

A significant advance in random matrix theory in recent years has been the development of a matrix transform theory based on spherical functions from harmonic analysis
for classes of random Hermitian matrices~\cite{KK2016,KR2016,KK2019,FKK2017}. One viewpoint of these studies is that they generalise to a matrix setting the approach to studying the distribution of sums and products of scalar random variables through the Fourier and Mellin transform respectively. Unitary invariance plays a key role, and an end product has been the identification of the previously unknown P\'olya ensembles --- intimately related to P\'olya frequency functions~\cite{Polya1913,Polya1915,Schoenberg} --- which exhibit a key closure property of the functional form of their joint eigenvalue probability density function (PDF) with respect to matrix addition or multiplication, as appropriate. The P\'olya ensembles are examples of determinantal point processes constructed out of a special class of biorthogonal functions. The latter permit explicit forms in terms of sums or integrals, which moreover allow for the correlation kernel to be written in a double contour integral form, which is a key ingredient in subsequent asymptotic analysis; see 
e.g.~\cite{LWZ16,FLT2020}. At a conceptual level, an explanation is thus provided not only for the persistence of determinantal structures of certain ensembles under matrix addition and multiplication, but also for aspects of its integrable structures.

The initial works considered products of unitary invariant positive definite Hermitian matrices~\cite{KK2016,KK2019,FKK2017}, and sums of unitary invariant Hermitian matrices
\cite{KR2016, Kieburg2019,FKK2017} from this viewpoint. Soon after works appeared involving antisymmetric matrices~\cite{FILZ2019,KFI2019} and Hermitian matrices with both positive and negative eigenvalues~\cite{Liu2017,FIL2018,Kieburg2020}.
In a recent work, 
three of the four authors of the present article introduced the spherical transform in the study of the randomised Horn problem for rank-1 additions and multiplications~\cite{ZKF2019}. This contained a class of unitary invariant matrix products of a type not studied in previous works using matrix transform methods, namely
multiplication on $\U(N)$. The harmonic analysis for this group goes back to Weyl and Harish-Chandra, as the spherical functions in the spherical transform~\cite{Helgason2000} are, up to a normalisation, the characters of the irreducible representations; those are the Schur polynomials in relation to $\U(N)$. In the present paper our aim is to identify the analogues of the P\'olya ensembles in this setting --- to be referred to as cyclic P\'olya ensemble --- and then to use the spherical transform 
to develop
theory culminating in the specification of the biorthogonal system and correlation kernel for the corresponding product matrices.

Definition~\ref{def:Polya-ensemble} of Section~\ref{sec:Polya}
specifies cyclic P\'olya ensembles as having PDF $p_N^{(U)}(z)$ of the eigenvalues
$z=\diag(z_1,\ldots,z_N)$ on the unit circle in the complex plane proportional to
\begin{equation}\label{Vandermonde}
\Delta_N(z) \det[(z_a\partial_a)^{b-1}\omega(z_a)]_{a,b=1\ldots,N}, \quad
 \Delta_N(z)=\det[z_a^{b-1}]_{a,b=1,\ldots,N}=\prod_{1\leq a< b\leq N}(z_b-z_a).
\end{equation}
Here the second expression for $\Delta_N(z)$ is the evaluation
of the Vandermonde determinant, and the weight function $\omega$ must belong to the class~\eqref{L1-set.b}. The latter
can be directly related to cyclic P\'olya frequency functions which have been defined for odd orders in~\cite{KRS1994}, and which we extend to even orders
in~Definition \ref{def:Polyafunction}. With $\mathbb I_N$ the set (\ref{Is}) and $s = (s_1,\dots, s_N) \in \mathbb I_N$, the spherical transform $\mathcal{S}^{(U)}(s)$
is specified in
Definition~\ref{spherical-transform}. From this, with $U$ drawn from a cyclic P\'olya ensemble with weight $\omega$, 
one calculates spherical transform
$\mathcal{S}\omega(s)$ (in this scalar case the spherical transform corresponds to the Fourier transform on an interval)
(c.f. Corollary~\ref{cor:spherical-Polya})
\begin{equation}\label{spherical-Polya}
		\mathcal{S}^{(U)}(s)=\prod_{j=1}^N\frac{\mathcal{S}\omega(s_j)}{\mathcal{S}\omega(j-1)}.
		\end{equation}
		
		With the aid of the spherical transform we are able to prove that the cyclic P\'olya ensembles
are closed under matrix multiplication. Specifically, with the weights of matrices from two cyclic P\'olya ensembles
being $\omega$ and $\hat{\omega}$, their product is a cyclic P\'olya ensemble with weight (part (2) of Theorem~\ref{thm:jpdfs})
\begin{equation}\label{new-omega}
		\widetilde{\omega}(z')=\widehat{\omega}\ast \omega(z')=\int_{\mathbb{S}_1}\frac{d\widetilde{z}}{2\pi \widetilde{z}} \widehat{\omega}\left(\frac{z'}{\widetilde{z}}\right)\omega(\widetilde{z})\in \widetilde{L}_N^1(\mathbb{S}_1)
		\end{equation}
		for all $z'\in\mathbb{S}_1$, the complex unit circle. Here the set $\widetilde{L}_N^1(\mathbb{S}_1)$ is specified by (\ref{L1-set.b}).
As another application of the spherical transform, a uniqueness theorem for the weight $\omega$
under the assumption of the non-vanishing of the first $N$ Laurent coefficients, counting from zero,
 is obtained.
 
 \begin{theorem}[Uniqueness of the Laurent Series and Weight]\label{thm:unique-Laurent}\
	
	Consider two cyclic P\'olya ensembles on $\U(N)$ associated to two weights $\omega_1$ and $\omega_2$ which have non-vanishing Laurent coefficients $u_{s_0}^{(1)},u_{s_0}^{(2)}\neq0$ for an integer $s_0\neq0,\ldots,N-1$. When their corresponding joint probability densities (see Eq.~\eqref{jpdf-Polya}) agree, the two weights can maximally differ by a global normalisation constant $C$. In particular, for $N=2M+1$ odd there is a $C>0$ and for $N=2M$ even there is a real $C\neq0$ with $\omega_1(z')=C\omega_2(z')$ for almost all $z'\in\mathbb{S}_1$.
\end{theorem}

The proof is in subsection~\ref{sec:Polya-def}. This result comes as quite a surprise in two ways, namely, firstly, the Haar measure weight is not unique for $N>1$
(Proposition~\ref{prop:cuep1}) while, secondly, any other cyclic P\'olya ensemble on $\U(N)$ is unique when normalising the weight. 

Part of the richness of the theory of products of unitary invariant positive definite Hermitian matrices is its tie in special functions by way of
the Meijer G-function~\cite{AI2015,KK2016,KK2019}. In the present setting, the role of the Meijer G-function is played by the bilateral
hypergeometric series defined as~\cite[Eq.~(16.4.16)]{NIST}
\begin{align}\label{biletral-hyper}
\phq\left[\left.
\begin{aligned}
a_1,\ldots,a_p\\
b_1,\ldots,b_q
\end{aligned}\right|x
\right]=\frac{\prod_{j=1}^q\Gamma(b_j)}{\prod_{j=1}^p\Gamma(a_j)}\sum_{s=-\infty}^\infty
\frac{\prod_{j=1}^p\Gamma(a_j+s)}{\prod_{j=1}^q\Gamma(b_j+s)}x^s.
\end{align}
The function $\phq$ is defined for all values of the variable $x$ such that $|x|=1$. If $x=-1$, we require $
\text{Re}(b_1+\cdots+b_q-a_1-\cdots-a_p)>1$, 
and if $x=1$, we require $\text{Re}(b_1+\cdots+b_q-a_1-\cdots-a_p)>0$.
Moreover, if any of the $a$ parameters is a negative integer or any of the $b$ parameters is a positive integer, then the series terminates above or below, respectively. If any of the $a$ parameters is a positive integer or if any of the $b$ parameters is a non-positive integer, the series is not defined as it experiences a pole.
The relevance of this class of special functions is immediate from their closure under multiplicative convolution on the complex unit circle
\begin{align}\label{HH}
\left[{_p}H_{q}\left[\left.
\begin{array}{c}
a_1,\ldots,a_p\\
b_1,\ldots,b_q
\end{array}\right|.\right]\ast{_{p'}}H_{q'}\left[\left.
\begin{array}{c}
c_1,\ldots,c_{p'}\\
d_1,\ldots,d_{q'}
\end{array}\right|.\right]\right](z')={_{p+p'}}H_{q+q'}\left[\left.
\begin{array}{c}
a_1,\ldots,a_p,c_1,\ldots,c_{p'}\\
b_1,\ldots,b_q,d_1,\ldots,d_{q'}
\end{array}\right|z\right].
\end{align}

Our primary example in this class is the cyclic P\'olya ensemble with weight
\begin{equation}\label{Jacobi-weight}
\omega_N^{\rm(Jac)}(z';\alpha,\gamma)=|(1+z')^{\alpha-2i\gamma}|(1+{z'}^*)^{N-1}={z'}^{-\alpha/2-i\gamma-N+1}(1+z')^{\alpha+N-1}, \quad
(\alpha > -1, \: \gamma \in \mathbb R)
\end{equation}
with $z'^*$ being the complex conjugate of $z'$. By applying partial differentiations to the Jacobi weight $\omega_N^{\rm(Jac)}$, equation (\ref{Vandermonde}) becomes proportional to 
\begin{equation}\label{CJi}
\frac{|\Delta_N(z)|^2}{\widetilde{C}_N}\prod_{j=1}^N\left|(1+z_j)^{\alpha-2i\gamma}\right|.
\end{equation}
This ensemble has also been considered in~\cite{FLT2020}, as an analogue of the Jacobi ensemble on the unitary group. It will be discussed further in subsection~\ref{sec:Jacobi}.

In Sec.~\ref{sec:prod} we compute the spectral statistics at finite matrix dimensions along the same lines as in~\cite{Kieburg2019,Kieburg2020}. In particular, we construct a bi-orthonormal pair of functions $\{ (P_j, Q_j) \}_{j=0,\dots,N-1}$ with which we can build the kernel of the corresponding determinantal point processes. 
The latter are given in series and integral forms. Specifically, for a cyclic polynomial ensemble an explicit series form of the bi-orthonormal pair of functions is given
in Proposition~\ref{prop:cyc.Pol}. Subject to a minor technical requirement on $\omega$, Corollary~\ref{cor:Christ-Polya} gives for the correlation kernel the integral form
\begin{equation}\label{kernel-cyc.Pol.c}
	K_N(z_1,z_2)=i\,\int_{0}^{2\pi}\frac{d\varphi}{2\pi}\varphi\, P_{N-1}(z_1e^{i \varphi})Q_{N}(z_2e^{i \varphi})+\frac{1-(z_1z_2^{-1})^{N}}{1-z_1z_2^{-1}}.
	\end{equation}
Here the integral on the right hand side has a form analogous to that known in the study of products of unitary invariant positive definite
Hermitian matrices~\cite{KK2016}, and recently shown to be key in studying the large $N$ hard edge asymptotics~\cite{FLT2020}. The second term is the kernel of the Circular Unitary Ensemble (CUE) which is the set of the unitary matrices distributed uniformly by the normalised Haar measure.

\section{Cyclic P\'olya Ensembles on $\U(N)$}\label{sec:Polya}

In subsection~\ref{sec:Haar-Polya}, we introduce the notion of cyclic polynomial ensembles, which are a natural generalisation of those on the real line~\cite{KZ2014}. Those ensembles exhibit the integrable structure of a determinantal point process~\cite{Borodin} and at the same time one set of functions of the corresponding bi-orthonormal pair are still polynomials. The utility of the latter is seen when computing the spectral transform of these ensembles, in subsection~\ref{sec:spherical}. As the polynomial part of the joint probability density is encoded in terms of a Vandermonde determinant, it cancels with the one from the spherical functions, which are in the present case the Schur polynomials. The problem is, however, that the product of two cyclic polynomial ensembles is not necessarily a cyclic polynomial ensemble again. The subclass which is closed under multiplicative matrix convolution are the cyclic P\'olya ensembles, introduced in subsection~\ref{sec:Polya-def}. As a benefit, those ensembles satisfy Harish-Chandra-like group integrals and have a closed multiplicative action on the set of cyclic polynomial ensembles for which we compute the resulting joint probability density of the eigenvalues. In the same subsection, we also prove that the weight of a cyclic P\'olya ensemble is unique if and only if it is not the Haar measure.

In subsection~\ref{sec:examples} we give several examples of cyclic P\'olya ensembles.
Therein we also show that we can readily construct cyclic P\'olya ensemble via products of certain exponentiated rank-1 random matrices. The class of P\'olya ensembles obtained in this way is by far exhaustive as can be seen by the circular Jacobi ensembles~\cite{WF2000,BD2002,FLT2020} for certain parameters.

The positivity condition of cyclic P\'olya ensembles is investigated in subsection~\ref{sec:CPF}. For this purpose, we extend the definition of cyclic P\'olya frequency functions on the circle~\cite{KRS1994} from odd to even orders. This is very important as we will see there is a subtle difference between these two kinds of dimensions which originates from the Vandermonde determinant.

\subsection{From the Haar Measure to Cyclic Polynomial Ensembles}\label{sec:Haar-Polya}

As stated in the Introduction, our aim is to advance the ideas of P\'olya ensembles~\cite{KR2016,KK2016,KK2019,FKK2017} for the additive and multiplicative matrix convolutions on spaces like the Hermitian matrices $\mathrm{Herm}(N)$ and the complex general linear group ${\rm GL}_{\mathbb{C}}(N)$ to the multiplicative convolution on the unitary matrices $\U(N)$. As we have learned from~\cite{KR2016,KK2016,KK2019,FKK2017}, those ensembles preserve the structure of determinantal point processes~\cite{Borodin} for their eigenvalue correlations under their respective matrix convolutions, i.e., the $k$-point correlation function has the form
\begin{equation}\label{k-point}
\begin{split}
R_k(z_1,\ldots,z_k)=&\frac{N!}{(N-k)!}\int_{\mathbb{S}_1^{N-k}}\frac{dz_{k+1}}{2\pi i z_{k+1}}\cdots \frac{dz_{N}}{2\pi i z_{N}}p_N(z)=\det[K_N(z_a,z_b)]_{a,b=1,\ldots, k}.
\end{split}
\end{equation}
The density $p_N(z)$ is the joint probability density of the eigenvalues $z=\diag(z_1,\ldots,z_N)$ on the torus $\mathbb{S}_1^N$ with $\mathbb{S}_1=\{z'\in\mathbb{C}|\, |z'|=1\}$ the centred complex unit circle and the reference measure $dz'/(2\pi i z')$, which is the normalised Haar measure on $\mathbb{S}_1$. We are interested in the unitarily invariant random matrix ensemble corresponding to $p_N(z)$ which is uniquely given because the Haar measure describing the distribution of the eigenvectors is unique. We recall that a function $f$ on $\U(N)$ is unitarily invariant if $f(U)=f(VUV^\dagger)$ for all $U,V\in\U(N)$ and $V^\dagger$ being the Hermitian adjoint of $V$.

The kernel $K_N(z_a,z_b)$ is, for instance, for the normalised Haar measure $d\mu(U)$ on $\U(N)$ of the form~\cite{Dyson}
\begin{equation}\label{kernel-Haar}
K_N^{\rm (Haar)}(z_a,z_b)=\sum_{j=0}^{N-1} \left(\frac{z_a}{z_b}\right)^j=\frac{1-(z_a/z_b)^N}{1-z_a/z_b}.
\end{equation}
The joint probability density of the eigenvalues of a Haar distributed unitary matrix is given by~\cite{Dyson}
\begin{equation}\label{jpdf-Haar}
p_N^{\rm (Haar)}(z)=\frac{1}{(2\pi)^N N!}|\Delta_N(z)|^2=\frac{(-1)^{N(N-1)/2}}{(2\pi)^N N!}\frac{\Delta_N^2(z)}{\det z^{N-1}},
\end{equation}
where $\det z:=\prod_{j=1}^N z_j$. The second equality is useful to identify~\eqref{jpdf-Haar} with the determinantal form of the Vandermonde determinant $\Delta_N(z)$ in~\eqref{Vandermonde}. Indeed the application of the generalised Andr\'eief indentity~\cite{KG2010,Andreief} on the second expression in~\eqref{jpdf-Haar} immediately yields the kernel~\eqref{kernel-Haar}. 
The latter reads that for suitably integrable sets of functions $\{P_{j-1}(z)\}_{j=1,\ldots,N}$ and $\{Q_{j-1}(z)\}_{j=1,\ldots,N}$, but otherwise arbitrary, we have
\begin{equation}\label{Andreief}
\begin{split}
&\int_{\mathbb{S}_1^{N-k}}\frac{dz_{k+1}}{2\pi i z_{k+1}}\cdots \frac{dz_{N}}{2\pi i z_{N}}\det[P_{b-1}(z_a)]_{a,b=1,\ldots,N}\det[Q_{b-1}(z_a)]_{a,b=1,\ldots,N}\\
& \qquad = (N-k)!\det\left[\begin{array}{cc} 0 & P_c(z_a) \\ -Q_d(z_b) & \displaystyle \int_{\mathbb{S}_1}\frac{dz'}{2\pi i z'} P_c(z')Q_d(z')\end{array}\right]_{\substack{a,b=1,\ldots k\\ c,d=0,\ldots,N-1}}.
\end{split}
\end{equation}
In the case of the Haar measure, one commonly chooses $P_c(z_a)=z_a^{c}$ and $Q_d(z_b)=z_b^{-d}$ to simplify the lower right block in the determinant on the right hand side of~\eqref{Andreief} to the identity. However the invariance of the determinant under linearly combining the rows and columns allows for a different basis. We will make use of this fact later on.

What we would like to concentrate on, now, is the generalisation of the joint probability density of the Haar measure~\eqref{jpdf-Haar} to a class of ensembles so that these densities satisfy the following conditions:
\begin{enumerate}
	\item	the joint probability density of the eigenvalues should have the form
	\begin{equation}\label{jpdf-poly}
	p_N(z)=\frac{1}{N!}\det[P_{b-1}(z_a)]_{a,b=1,\ldots,N}\det[Q_{b-1}(z_a)]_{a,b=1,\ldots,N}
	\end{equation}
	so that it is guaranteed that the eigenvalue statistics build a determinantal point process;
	\item	the span of $\{P_{j-1}(z)\}_{j=1,\ldots,N}$ is still the vector space of polynomials of order $N-1$;
	\item	when $U_1,U_2\in\U(N)$ are two independent, not necessarily identically distributed, unitarily invariant random matrices with joint probability densities of their eigenvalues of the form~\eqref{jpdf-poly}, then, also the eigenvalues of $U_1U_2$ are distributed along the form~\eqref{jpdf-poly}. Certainly, the functions $Q_{b-1}(z_a)$ may vary for $U_1$, $U_2$ and $U_1U_2$.
\end{enumerate}

The first two conditions bring us to our first definition of the notion of a polynomial ensemble on $\U(N)$ which is the counterpart of polynomial ensembles for real spectra~\cite{KZ2014}. For this purpose we define the set of functions
\begin{equation}\label{L1-set}
L_N^1(\mathbb{S}_1)=\{w\in L^1(\mathbb{S}_1)|\, [w(z)]^*=z^{N-1}w(z)\}.
\end{equation}
We note that the $L^1$-functions on the complex unit circle are all functions that are absolutely integrable with respect to the Haar measure $|dz/z|=d\vartheta$ on $\mathbb{S}_1$ with $z=e^{i\vartheta}$ and $\vartheta\in]-\pi,\pi[$. The necessity of the condition $[w(z)]^*=z^{N-1}w(z)$ results from the following definition.

\begin{definition}[Cyclic Polynomial Ensemble]\label{def:polynomial}\
	
	A unitarily invariant random matrix $U\in\U(N)$ is called a cyclic polynomial ensemble associated to the weights $\{w_j\}_{j=0,\ldots,N-1}\subset L_N^1(\mathbb{S}_1)$ iff its joint probability distribution of its eigenvalues $z=\diag(z_1,\ldots,z_N)\in\mathbb{S}_1^N$ has the form
	\begin{equation}\label{jpdf-polynomial}
	p_N^{(U)}(z)=\frac{1}{C_N N!}\frac{\Delta_N(z)}{i^{N(N-1)/2}}\det[w_{b-1}(z_a)]_{a,b=1,\ldots,N}\geq0
	\end{equation}
	with respect to the measure $\prod_{j=1}^N dz_j/(2\pi i z_j)$ and the normalisation constant
	\begin{equation}\label{norm-polynomial}
	C_N=\det\left[\int_{\mathbb{S}_1}\frac{dz'}{2\pi i z'} (-i z')^{a-1}w_{b-1}(z')\right]_{a,b=1,\ldots,N}>0.
	\end{equation}
\end{definition}

One can see that the additional condition in the set~\eqref{L1-set} guarantees that the joint probability density is real because of
\begin{equation}\label{Vand-ident}
[\Delta_N(z)]^*=(-1)^{N(N-1)/2}\frac{\Delta_N(z)}{\prod_{j=1}^N z_j^{N-1}},
\end{equation} 
which we have already exploited for the second equality in~\eqref{jpdf-Haar}. Hence, there is certainly also a real representation when choosing the coordinates $z_j=e^{i\theta_j}$ with $\theta_j\in]-\pi,\pi[$ of the form
\begin{equation}\label{jpdf-polynomial.b}
\begin{split}
&p_N^{(U)}(e^{i\theta})=\frac{1}{C_N N!}\left(\prod_{1\leq a<b\leq N}2\sin\left[\frac{\theta_a-\theta_b}{2}\right]\right)\det[\widehat{w}_{b-1}(\theta_a)]_{a,b=1,\ldots,N}\\
&{\rm with}\ \widehat{w}_{b-1}(\theta_a)=[\widehat{w}_{b-1}(\theta_a)]^*=e^{i(N-1)\theta_a/2}w_{b-1}(e^{i\theta_a}).
\end{split}
\end{equation}
For the square root of the complex phases, the branch cut is taken along the negative real axis. The price that we have to pay is that for even dimensions $N$ the functions $\widehat{w}_{b-1}(\theta_a)=\widehat{w}_{b-1}(\theta_a+4\pi)$ are only $4\pi$ periodic, more precisely they are $2\pi$ anti-periodic, $\widehat{w}_{b-1}(\theta_a)=-\widehat{w}_{b-1}(\theta_a+2\pi)$, not like the $2\pi$ periodicity for odd $N$. Indeed, the $2\pi$ periodicity is always preserved for the weights $w_{b-1}(e^{i\theta})=w_{b-1}(e^{i(\theta+2\pi)})$. Hence, this change of periodicity is not a problem, the joint probability density $p_N^{(U)}(e^{i\theta})$ stays always $2\pi$ periodic in each angle $\theta_j$. This observation has some important consequences in the explicit representation of some ensembles as it has been already noted in~\cite{Liechty}, and we will see this below, too. 

In contrast, the positivity of the joint probability density cannot be so easily traced back and ensured. We will discuss this in more detail for the P\'olya ensembles on $\U(N)$ that have to be still defined, yet.

\subsection{Spherical Transforms on $\U(N)$}\label{sec:spherical}

Let us turn our attention to the last of the three aforementioned conditions, namely that the product $U_1U_2$ of two independent, unitarily invariant random matrices $U_1,U_2\in\U(N)$, that are drawn from two (maybe different) cyclic polynomial ensembles, is also a cyclic polynomial ensemble. We will see in the ensuing discussion that this is not true for two arbitrary cyclic polynomial ensembles. We emphasize that the unitary invariance of the product is a direct consequence of the one of $U_1$ and $U_2$ because of $VU_1U_2V^\dagger=(VU_1V^\dagger)(VU_2V^\dagger)$.

The tool we need to discuss products of unitary matrices is the result of a successful combination of harmonic analysis and group and representation theory; it is the method of spherical transforms~\cite{Helgason2000}. For the multiplicative action on the unitary group $\U(N)$ this was recently introduced in RMT and applied to the multiplicative Horn problem by some of the present authors in~\cite{ZKF2019}. We will briefly repeat the definition of the spherical transform and recall some of its properties. To this aim we define the multi-index set
\begin{equation}\label{Is}
\mathbb{I}_N=\{(s_1,\ldots,s_N)\in\mathbb{Z}^N| s_a\neq s_b\ {\rm when}\ a\neq b\}.
\end{equation}

\begin{definition}[Spherical Transform]\
	
	Let $s=(s_1,\ldots,s_N)\in\mathbb{I}_N$. The spherical transform of an $L^1$-function $f(U)$ on $\U(N)$ is given by
	\begin{equation}\label{spherical-transform}
	\mathcal{S}f(s)=\int_{\U(N)} d\mu(U) f(U)\Phi(U;s),
	\end{equation}
	with $d\mu(U)$ the normalised Haar measure on $\U(N)$. The spherical function given by
	\begin{equation}\label{spherical-function}
	\Phi(U;s)=\frac{{\rm ch}_s(U)}{{\rm ch}_s(\mathbf{1}_N)}=\left(\prod_{j=0}^{N-1}j!\right)\frac{\det[z_a^{s_b}]_{a,b=1\ldots,N}}{\Delta_N(z)\Delta_N(s)}
	\end{equation}
	is the ratio of the character of $U$ and the $N$-dimensional identity $\mathbf{1}_N$, where $\mathrm{ch}_s$ denotes the character of $U$ which is the trace of an irreducible representation of $U$ (see e.g.~\cite[Ch. IV Sec. I]{Helgason2000}). The right hand side of~\eqref{spherical-function} is given in terms of the eigenvalues $z=\diag(z_1,\ldots,z_N)\in\mathbb{S}_1^N$ of the matrix $U$.
\end{definition}



For the additive and multiplicative convolution on ${\rm Herm}(N)$, ${\rm GL}_{\mathbb{C}}(N)$ etc., the points where two or more indices $s_j$ of $s=(s_1,\ldots,s_N)$ agree with each other are of measure zero, and therefore one may exclude those points when inverting the spherical transform. For the multiplication of unitary matrices the Fourier space of the ``frequencies'' $s$ has to be discrete since the unitary group is compact. The natural measure on $\mathbb{Z}^N$ is the Dirac measure. Hence, it is still crucial to exclude those points as they will be not of measure zero. The deeper representation theoretical reason is that the characters of finite dimensional irreducible representations of compact groups must be polynomials of all matrix entries of the group element $U$. In our case these are the Schur polynomials . The frequencies $s$ are related to the partition of the corresponding irreducible representation. When two $s_j$ agree we are forced to understand the character by l'H\^opital's rule creating logarithms of the eigenvalues $z_a$ which are not any more polynomials of the matrix $U$. Therefore, these terms must be excluded to agree with the group theoretical insights.

We would like to also point out that the function $f$ does not necessarily need to be unitarily invariant. However, when it is unitarily invariant the formula~\eqref{spherical-transform} immediately simplifies to
\begin{equation}
\mathcal{S}f(s)=\frac{\prod_{j=0}^{N-1}j!}{N!}\int_{\mathbb{S}_1^N} \left(\prod_{j=1}^N\frac{dz_j}{2\pi i z_j}\right)|\Delta_N(z)|^2 f(z)\frac{\det[z_a^{s_b}]_{a,b=1\ldots,N}}{\Delta_N(z)\Delta_N(s)}.
\end{equation}
With a slight abuse of notation we also write
\begin{equation}\label{spherical-transform-jpdf}
\mathcal{S}p_N^{(U)}(s)=\left(\prod_{j=0}^{N-1}j!\right)\int_{\mathbb{S}_1^N} \left(\prod_{j=1}^N\frac{dz_j}{2\pi i z_j}\right)p_N^{(U)}(z)\frac{\det[z_a^{s_b}]_{a,b=1\ldots,N}}{\Delta_N(z)\Delta_N(s)},
\end{equation}
where now $p_N^{(U)}(z)$ is the joint probability density of the eigenvalues which comprises a major part of the Haar measure on $\U(N)$ since its reference measure is the Haar measure $\prod_{j=1}^N dz_j/(2\pi i z_j)$ on the $N$-dimensional torus $\mathbb{S}_1^N$.

\begin{remark}[Probability Densities on $\U(N)$ and $\mathbb{S}_1^N$]\
	
	To distinguish the probability density of a unitarily invariant random matrix $U$ on $\U(N)$ with the joint probability density function of the eigenvalues on the torus $\mathbb{S}_1^N$ we apply the following notation.
	\begin{enumerate}
		\item	{\bf The probability density of $U\in\U(N)$} is denoted by $f_N^{(U)}$ where the superscript indicates the random matrix it corresponds to. The reference measure is the normalised Haar measure $d\mu(U')$ on $\U(N)$. In particular the density is normalised as follows
		\begin{equation}
		\int_{\U(N)}d\mu(U') f_N^{(U)}(U')=1.
		\end{equation}
		Therefore, the Haar measure on $\U(N)$ has the probability density $f_N^{\rm(Haar)}(U')=1$.
		\item	{\bf The joint probability density of the eigenvalues} $z=\diag(z_1,\ldots,z_N)\in\mathbb{S}_1^N$ of the random matrix $U$ is coined $p_N^{(U)}$ and is normalised with respect to the normalised Haar measure on $\mathbb{S}_1^N$, i.e.,
		\begin{equation}
		\int_{\mathbb{S}_1^N} \left(\prod_{j=1}^N\frac{dz_j}{2\pi i z_j}\right) p_N^{(U)}(z)=1.
		\end{equation}	
		For the Haar measure the corresponding joint probability density of the eigenvalues is given in~\eqref{jpdf-Haar}.
		\item	{\bf The relation} between a unitarily invariant density $f_N^{(U)}$ and $p_N^{(U)}$ is given by
		\begin{equation}
		p_N^{(U)}(z)=\frac{1}{N!}|\Delta_N(z)|^2f_N^{(U)}(z).
		\end{equation}
		Therefore, the spherical transform of $f_N^{(U)}$ agrees with the one of $p_N^{(U)}$,
		\begin{equation}\label{spherical-relation}
		\mathcal{S}f_N^{(U)}=\mathcal{S}p_N^{(U)}=\mathcal{S}^{(U)}.
		\end{equation}
		The abbreviation $\mathcal{S}^{(U)}$ highlights this feature. We make use of it when we do not need to highlight which density we consider.
	\end{enumerate}
\end{remark}

As the spherical transform plays a crucial role in the ensuing sections, we would like to summarise some of its properties, see~\cite{Helgason2000,ZKF2019}.
\begin{enumerate}
	\item	{\bf The normalisation} is given by $s=s^{(0)}$ with $s_j^{(0)}=j-1$ because of
	\begin{equation}
	\Phi(U;s^{(0)})=1
	\end{equation}
	so that we have
	\begin{equation}\label{spherical-normalisation}
	\mathcal{S}f(s^{(0)})=\int_{\U(N)} d\mu(U) f(U),
	\end{equation}
	which equals $1$ when $f$ is a probability density on $\U(N)$.
	\item	{\bf The inverse} of the spherical transform is for unitarily invariant ensembles guaranteed when restricting to the image of $\mathcal{S}$ and it is explicitly given by~\cite[Proposition 1 in Sec.~4.2]{ZKF2019}
	\begin{equation}\label{spherical-inverse}
	\begin{split}
	\mathcal{S}^{-1}[\mathcal{S}f](U)=&\frac{1}{N!\prod_{j=0}^{N-1}(j!)^2}\lim_{t\to0^+}\sum_{s\in\mathbb{I}_N}\Delta_N^2(s) \mathcal{S}f(s)\Phi(U^\dagger;s)\\
	&\times\exp\left[-t\tr \left(s+\frac{1-N}{2}\mathbf{1}_N\right)^2+t\sum_{j=0}^{N-1} \left(j+\frac{1-N}{2}\right)^2\right].
	\end{split}
	\end{equation}
	The regularisation $\exp\left[-t\tr \left(s+\frac{1-N}{2}\mathbf{1}_N\right)^2+t\sum_{j=0}^{N-1} \left(j+\frac{1-N}{2}\right)^2\right]$ is only important for those $L^1$-functions for which the series of $\Delta_N^2(s) \mathcal{S}f(s)$ on $s\in\mathbb{I}_N$ is not absolutely convergent. In cases where the absolute convergence is given, we can neglect this auxiliary term. We would like to underline that $\mathcal{S}^{-1}[\mathcal{S}f](U)$ and $f(U)$ only need to agree almost everywhere as it is known that there might be inconsistencies at points where $f$ is discontinuous. Those points, however, are irrelevant when the reference measure is the Haar measure on $\U(N)$.
	\item	The spherical transform is evidently {\bf symmetric} in its arguments $s$ because of the symmetry of the spherical function $\Phi(U;s)=\Phi(U;s_\pi)$ for any permutation $s_\pi$ of the multi-index $s\in\mathbb{I}_N$.
	\item	{\bf The factorisation theorem} makes statements on the spherical transform of the random matrix $U_1U_2$ where $U_1\in\U(N)$ is fixed and $U_2\in\U(N)$ is a unitarily invariant random matrix. Say $f_N^{(U_2)}$ and $f_N^{(U_1U_2)}$ are the respective probability densities on $\U(N)$. Then, we have
	\begin{equation}\label{convolution-factorisation-fixed}
	\mathcal{S}^{(U_1U_2)}(s)=\Phi(U_1;s)\mathcal{S}^{(U_2)}(s).
	\end{equation}
	This equation also holds when $U_1=V\tilde{U}_1V^\dagger$ with $\tilde{U}_1\in\U(N)$ fixed and $V\in\U(N)$ Haar distributed because characters and, hence, the spherical function are invariant under cyclic permutations, i.e., ${\rm ch}_s(AB)={\rm ch}_s(BA)$; it is a trace of the product of $A$ and $B$ in a certain irreducible representation of the unitary group. Thence, $V\tilde{U}_1V^\dagger U_2$ and $\tilde{U}_1V^\dagger U_2V$ and, therefore, $\tilde{U}_1 U_2$ (because of the unitarily invariance of $U_2$) share the same joint probability density of the eigenvalues.
	
	Equation~\eqref{convolution-factorisation-fixed} is a direct consequence for the well-known factorisation formula for characters,
	\begin{equation}
	\int_{\U(N)}d\mu(U){\rm ch}_s(U_1UU_2U^\dagger)=\frac{{\rm ch}_s(U_1){\rm ch}_s(U_2)}{{\rm ch}_s(\mathbf{1}_N)}.
	\end{equation}
	When also the matrix $U_1$ is a random matrix on $\U(N)$ drawn from the probability density $f_N^{(U_1)}$, Eq.~\eqref{convolution-factorisation-fixed} reads then
	\begin{equation}\label{convolution-factorisation-random}
	\mathcal{S}^{(U_1U_2)}(s)=\mathcal{S}^{(U_1)}(s)\mathcal{S}^{(U_2)}(s).
	\end{equation}
	The multiplicative convolution on $\U(N)$,
	\begin{equation}\label{mult-conv}
	\begin{split}
	f_N^{(U_1U_2)}(U)=f_N^{(U_1)}\ast f_N^{(U_2)}(U)&=\\
	\int_{\U(N)}d\mu(U')f_N^{(U_1)}(U')& f_N^{(U_2)}(U{U'}^\dagger)=f_N^{(U_2)}\ast f_N^{(U_1)}(U),
	\end{split}
	\end{equation}
	can be also rewritten into form
	\begin{equation}\label{convolution-spherical}
	f_N^{(U_1U_2)}(U)=f_N^{(U_1)}\ast f_N^{(U_2)}(U)=\mathcal{S}^{-1}\left[\mathcal{S}f_N^{(U_1)}\mathcal{S}f_N^{(U_2)}\right](U).
	\end{equation}
	This is one effective way to evaluate a convolution and of which we will rely later on. We would like to point out that when $N=1$ equation~\eqref{convolution-factorisation-random} reduces to
	\begin{equation}\label{U(1)_convolve}
	\mathcal{S}w_{1}(s')\mathcal{S}w_2(s')=\mathcal{S}[w_{1}\ast w_2](s')
	\end{equation}
	for $w_1,w_2$ two probability densities defined on $\mathbb S_1$ and $s'\in\mathbb{Z}$. 
\end{enumerate}

\begin{remark}\
	
	Certainly, the relations above also carry over to the spherical transform of the joint probability density of the eigenvalues $z=\diag(z_1,\ldots,z_N)$ of $U\in\U(N)$, due to~\eqref{spherical-relation}. Especially, the inverse of the spherical transform is then explicitly~\cite[Lemma 3 in Sec.~4.2]{ZKF2019}
	\begin{equation}\label{spherical-inverse-jpdf}
	\begin{split}
	\mathcal{S}^{-1}[\mathcal{S}p_N^{(U)}](z)=&\frac{|\Delta_N(z)|^2}{N!\prod_{j=0}^{N-1}(j!)^2}\lim_{t\to0^+}\sum_{s\in\mathbb{I}_N}\Delta^2(s) \mathcal{S}p_N^{(U)}(s)\Phi(z^{-1};s)\\
	&\times\exp\left[-t\tr \left(s+\frac{1-N}{2}\mathbf{1}_N\right)^2+t\sum_{j=0}^{N-1} \left(j+\frac{1-N}{2}\right)^2\right].
	\end{split}
	\end{equation}
	We will mostly work on the level of the eigenvalues, in the ensuing sections, so that Eqs.~\eqref{spherical-transform-jpdf} and~\eqref{spherical-inverse-jpdf} will be of importance for us.
\end{remark}

As a simple exercise we will first compute the spherical transform of an arbitrary cyclic polynomial ensemble.

\begin{proposition}[Spherical Transform of a Cyclic Polynomial Ensemble]\label{prop:spher.polynomial}\
	
	The spherical transform of the cyclic polynomial ensemble in Definition~\ref{def:polynomial} with $p_N^{(U)}(z)$ the joint probability density~\eqref{jpdf-polynomial} of the eigenvalues $z$ is given by
	\begin{equation}\label{spherical-polynomial}
	\mathcal{S}^{(U)}(s)=\mathcal{S}p_N^{(U)}(s)=\frac{\prod_{j=0}^{N-1}j!}{\Delta_N(s)}\frac{\det[\mathcal{S}w_{b-1}(s_a)]_{a,b=1,\ldots N}}{\det[\mathcal{S}w_{b-1}(a-1)]_{a,b=1,\ldots N}}
	\end{equation}
	for all $s=\diag(s_1,\ldots,s_N)\in\mathbb{I}_N$. The spherical transform for the weights is given by the univariate Fourier transform
	\begin{equation}\label{univariate-spherical}
	\mathcal{S}w_{b-1}(s_a)=\int_{\mathbb{S}_1}\frac{dz'}{2\pi i z'}{z'}^{s_a}w_{b-1}(z').
	\end{equation}
\end{proposition}

Due to the invertibility of the spherical transform one can give a stronger statement and say that a unitarily invariant random matrix $U\in\U(N)$ is drawn from a cyclic polynomial ensemble iff its spherical transform has the form~\eqref{spherical-polynomial}. Here one needs to restrict the domain of the $\mathcal{S}^{-1}$ to the image of $\mathcal{S}$ for unitarily invariant probability densities on $\U(N)$ with respect to the Haar measure.

{\bf Proof of Proposition~\ref{prop:spher.polynomial}:}

The constant $C_N$, see~\eqref{norm-polynomial}, obviously accounts for the denominator in~\eqref{spherical-polynomial} when employing the definition~\eqref{univariate-spherical} of the univariate Fourier transform. Thus, we get
\begin{equation}
\mathcal{S}^{(U)}(s)=\frac{\prod_{j=0}^{N-1}j!}{N!\det[\mathcal{S}w_{b-1}(a-1)]_{a,b=1,\ldots N}}\int_{\mathbb{S}_1^N} \left(\prod_{j=1}^N\frac{dz_j}{2\pi i z_j}\right)\det[w_{b-1}(z_a)]_{a,b=1,\ldots,N}\frac{\det[z_a^{s_b}]_{a,b=1\ldots,N}}{\Delta_N(s)}
\end{equation}
after cancelling some phase factors and the Vandermonde determinants $\Delta_N(z)$. Applying the original Andr\'eief identity~\cite{Andreief} which is~\eqref{Andreief} for $k=0$ and employing anew Eq.~\eqref{univariate-spherical} we arrive at~\eqref{spherical-polynomial}.
\hfill$\square$

As a trivial consequence, we obtain the following corollary for the Haar measure. Using the Andr\'eief identity, one only needs to identify $w_{b-1}(z_a)=z_a^{1-b}$ and carry out the integral which yields Kronecker deltas of the form $\delta_{s_a,b-1}$. The determinant tells us that $(s_1,\ldots,s_N)$ has to be a permutation of $(0,\ldots,N-1)$. Therefore, the constant and the sign in~\eqref{spherical-polynomial} cancel each other.

\begin{corollary}[Spherical Transform of the Haar Measure]\label{cor:spher.Haar}\
	
	The spherical transform of the Haar measure is
	\begin{equation}\label{spherical-Haar}
	\mathcal{S}p_N^{\rm (Haar)}(s)=\prod_{j=1}^N\chi_{[0,N-1]}(s_j)
	\end{equation}
	for all $s=\diag(s_1,\ldots,s_N)\in\mathbb{I}_N$, where $\chi_{[0,N-1]}(s_j)$ is the indicator function on the interval $[0,N-1]$, meaning it is only $1$ when $s_j\in[0,N-1]$ and vanishes otherwise.
\end{corollary}

\subsection{Cyclic P\'olya Ensembles}\label{sec:Polya-def}

Considering two random matrices $U_1,U_2\in\U(N)$ drawn from the probability densities $f_N^{(U_1)}$ and $f_N^{(U_2)}$, we readily notice that their product $U_1U_2$ do not necessarily yield a cyclic polynomial ensemble even if they were both cyclic polynomial ensembles. Say $U_1$ is associated to the weights $\{w_j^{(1)}\}_{j=0,\ldots,N-1}\subset L_N^1(\mathbb{S}_1)$ and $U_2$ is associated to $\{w_j^{(2)}\}_{j=0,\ldots,N-1}\subset L_N^1(\mathbb{S}_1)$. Then, the spherical transform of the probability density $f_N^{(U_1U_2)}$ for the product $U_1U_2$ is given by
\begin{equation}
\mathcal{S}^{(U_1U_2)}(s)=\frac{\prod_{j=0}^{N-1}(j!)^2}{\Delta_N^2(s)}\frac{\det[\mathcal{S}w_{b-1}^{(1)}(s_a)]_{a,b=1,\ldots N}}{\det[\mathcal{S}w_{b-1}^{(1)}(a-1)]_{a,b=1,\ldots N}}\frac{\det[\mathcal{S}w_{b-1}^{(2)}(s_a)]_{a,b=1,\ldots N}}{\det[\mathcal{S}w_{b-1}^{(2)}(a-1)]_{a,b=1,\ldots N}}.
\end{equation}
The weights have to satisfy certain conditions so that this product simplifies to the form~\eqref{spherical-polynomial}. The simplest way to reach this goal is that one of the two determinants, say the one for $U_2$ in the numerator can be reduced to the form
\begin{equation}
\det[\mathcal{S}w_{b-1}^{(2)}(s_a)]_{a,b=1,\ldots N}=\Delta_N(s)\prod_{j=1}^N\sigma(s_j)
\end{equation}
with $\sigma$ being a complex valued function on $\mathbb{Z}$. Note that the symmetries in the argument $s$ need to be preserved for the ansatz which is here the case. Without loss of generality, one can say that we have
\begin{equation}
\mathcal{S}w_{b-1}^{(2)}(s_a)=q_{b-1}(s_a)\sigma(s_a),
\end{equation}
with $q_{b-1}(s_a)=s_a^{b-1}+\ldots$ a monic polynomial of order $b-1$ or when applying the inverse spherical transform, we arrive at
\begin{equation}
w_{b-1}^{(2)}(z')=\mathcal{S}^{-1}[q_{b-1}(s_a)\sigma(s_a)](z')=q_{b-1}(-z'\partial_{z'})\mathcal{S}^{-1}\sigma(z').
\end{equation}
Here, we used the identity
\begin{equation}\label{derivative-spherical}
\mathcal{S}[-z'\partial_{z'}f(z')](s')=s'\mathcal{S}f(s')
\end{equation}
for any suitably differentiable and integrable function $f$ on $\mathbb{S}_1$. It is a direct consequence of~\eqref{univariate-spherical}. From this perspective it is very natural to define a subclass of cyclic polynomial ensembles on $\U(N)$, namely cyclic P\'olya ensembles. Their name is born out from their relation to P\'olya frequency functions on the complex unit circle which will be discussed in subsection~\ref{sec:CPF}. For this aim, we need to define the functions
\begin{equation}\label{L1-set.b}
\widetilde{L}_N^1(\mathbb{S}_1)=\{w\in L_N^1(\mathbb{S}_1)|\, w\ \text{is $(N-1)$-times differentiable,}\ \partial^jw\in L^1(\mathbb{S}_1)\ \text{for all}\ j=0,\ldots,N-1\}.
\end{equation}
Let us highlight that the functions are only $(N-2)$-times continuous differentiable while its $N-1$ needs only to exist almost everywhere.

\begin{definition}[Cyclic P\'olya Ensemble]\label{def:Polya-ensemble}\
	
	A unitarily invariant random matrix $U\in\U(N)$ is drawn from a cyclic P\'olya ensemble on $\U(N)$ associated to the weight $\omega\in \widetilde{L}_N^1(\mathbb{S}_1)$ iff its joint probability density of its eigenvalues $z=\diag(z_1,\ldots,z_N)\in\mathbb{S}_1^N$ can be written in the form
	\begin{equation}\label{jpdf-Polya}
	p_N^{(U)}(z)=\frac{1}{ N! \prod_{j=0}^{N-1}[j!\mathcal{S}\omega(j)]}\Delta_N(z)\det[(-z_a\partial_a)^{b-1}\omega(z_a)]_{a,b=1,\ldots,N}\geq0.
	\end{equation}
	Hereafter, $\partial_a$ is an abbreviation for $\partial_{z_a}$.
\end{definition}

One can readily check the normalisation and that $ p_N^{(U)}$ is real-valued. For instance, the Andr\'eief integral identity~\cite{Andreief}, see~\eqref{Andreief} for $k=0$, leads to
\begin{equation}
\begin{split}
\int_{\mathbb{S}_1^N} \left(\prod_{j=1}^N\frac{dz_j}{2\pi i z_j}\right) p_N^{(U)}(z)=&\frac{\det[\int_{\mathbb{S}_1}dz'/(2\pi i z')\ {z'}^{a-1} (-z'\partial)^{b-1}\omega(z')]_{a,b=1,\ldots,N}}{\prod_{j=0}^{N-1}[j!\mathcal{S}\omega(j)]}\\
=&\frac{\det[(a-1)^{b-1}\mathcal{S}\omega(a-1)]_{a,b=1,\ldots,N}}{\prod_{j=0}^{N-1}[j!\mathcal{S}\omega(j)]}\\
=&1.
\end{split}
\end{equation}
The realness results from $[\omega(z_a)]^*=z_a^{N-1}\omega(z_a)$ and $(z_a\partial_a)^*=-z_a\partial_a$ for all $z_a\in\mathbb{S}_1$. The minus sign cancels with the minus sign in~\eqref{Vand-ident} and the factors of $z_a^{N-1}$ come in handy too. Thus we note that the commutation of the factor ${z_a}^{N-1}$ with $(z_a\partial_a)^j$ yields a monic polynomial in $(z_a\partial_a)^j$ of order $j$, i.e.,
\begin{equation}
(z_a\partial_a)^jz_a^{N-1}=z_a^{N-1}(z_a\partial_a+N-1)^j.
\end{equation}
After a linear combination of the rows in the determinant we arrive at the same determinant again.

As a simple consequence of the definition of a P\'olya ensemble and Proposition~\ref{prop:spher.polynomial} the spherical transform can be made
explicit. One only needs to replace $w_{b-1}$ by $(-z'\partial)^{b-1}\omega$ in~\eqref{spherical-polynomial} and to exploit~\eqref{derivative-spherical}.

\begin{corollary}[Spherical Transform]\label{cor:spherical-Polya}\
	
	\begin{enumerate}
		\item	The spherical transform of the cyclic P\'olya ensemble in Definition~\ref{def:Polya-ensemble} is equal to (\ref{spherical-Polya}).
		for all $s=\diag(s_1,\ldots,s_N)\in\mathbb{I}_N$..
		\item	The spherical transform of the inverse random matrix $U^{-1}\in\U(N)$ of part (a) is
		\begin{equation}\label{spherical-Polya-inv}
		\mathcal{S}^{(U^{-1})}(s)=\prod_{j=1}^N\frac{\mathcal{S}\omega(N-s_j-1)}{\mathcal{S}\omega(N-j)}
		\end{equation}
		for all $s=\diag(s_1,\ldots,s_N)\in\mathbb{I}_N$.
		Therefore, $U^{-1}$ is drawn from a P\'olya ensemble, too, with the weight
		\begin{equation}\label{weight-inverse}
		\widetilde{\omega}(z')= z'^{1-N}\omega({z'}^{-1})=[\omega({z'}^{-1})]^*.
		\end{equation}
	\end{enumerate}
\end{corollary}

{\bf Proof of Corollary~\ref{cor:spherical-Polya}:}

As already mentioned, Eq.~\eqref{spherical-Polya} is a very direct consequence of Eq.~\eqref{spherical-polynomial}.
The second statement, in contrast, follows from the fact that the inverse of a unitary matrix implies that we consider the inverse of its eigenvalues such that their joint probability density is given by replacing $z\leftrightarrow z^{-1}$ in the original joint probability density~\eqref{jpdf-Polya}. This immediately leads to~\eqref{weight-inverse} and, hence, Eq.~\eqref{spherical-Polya}.
\hfill$\square$

With the aid of this result we come back to products involving cyclic P\'olya ensembles which has been the motivation from the start and has led us to the introduction of this class of unitary random matrices. The following Theorem is our first main result.

\begin{theorem}[Products involving P\'olya Ensembles]\label{thm:jpdfs}\
	
	Let $U_1$ be a unitarily invariant random matrix drawn from a cyclic P\'olya ensemble on $\U(N)$, associated to the weight $\omega\in \widetilde{L}_N^1(\mathbb{S}_1)$.
	\begin{enumerate}
		\item	Drawing a second unitarily invariant random matrix $U_2\in\U(N)$ from a cyclic polynomial ensemble associated to the weights $\{w_{j}\}_{j=0,\ldots,N-1}\subset L_N^1(\mathbb{S}_1)$. Then, $U=U_1U_2$ belongs to a cyclic polynomial ensemble associated to the weights
		\begin{equation}\label{new-w}
		\widetilde{w}_{j}(z')=w_j\ast \omega(z')=\int_{\mathbb{S}_1}\frac{d\widetilde{z}}{2\pi \widetilde{z}} w_j\left(\frac{z'}{\widetilde{z}}\right)\omega(\widetilde{z})\in L_N^1(\mathbb{S}_1)
		\end{equation}
		for all $j=0,\ldots,N-1$ and $z'\in\mathbb{S}_1$.
		\item	Choosing a second unitarily invariant $U_2\in\U(N)$ from a cyclic P\'olya ensemble associated to the weight $\widehat{\omega}\in \widetilde{L}_N^1(\mathbb{S}_1)$, $U=U_1U_2$ is a cyclic P\'olya ensemble associated to the weight (\ref{new-omega}).
		\item	Let $U_2\in\U(N)$ be fixed with the pair-wise different eigenvalues $x=\diag(x_1,\ldots,x_N)\in\mathbb{S}_1^N$, i.e., $x_a\neq x_b$ for $a\neq b$ and $V\in\U(N)$ should be Haar distributed. Then, the random matrix $U=U_1VU_2V^\dagger$ belongs to a cyclic polynomial ensemble associated to the weights
		\begin{equation}
		\widetilde{w}_{j}(z')=\omega\left(\frac{z'}{x_{j+1}}\right)
		\end{equation}
		for all $j=0,\ldots,N-1$ and $z'\in\mathbb{S}_1$. In particular, the joint probability density of the eigenvalues $z=\diag(z_1,\ldots,z_N)\in\mathbb{S}_1^N$ of $U$ is equal to
		\begin{equation}\label{fixed.jpdf}
		p_N^{(U)}(z|x)=\frac{1}{N!\prod_{j=0}^{N-1}\mathcal{S}\omega(j)}\frac{\Delta_N(z)}{\Delta_N(x)}\det\left[\omega\left(\frac{z_a}{x_{b}}\right)\right]_{a,b=1\ldots,N}.
		\end{equation}
		For a degenerate spectrum of $U_2$, one needs to apply l'H\^opital's rule.
	\end{enumerate}
\end{theorem}

{\bf Proof of Theorem~\ref{thm:jpdfs}:}

The first two statements are straightforward consequences of the bijectivity of the spherical transform and the factorisation identity~\eqref{convolution-factorisation-random}. Explicitly, the spherical transform of $U=U_1U_2$ is
\begin{equation}
\mathcal{S}^{(U)}(s)=\frac{\prod_{j=0}^{N-1}j!}{\Delta_N(s)}\frac{\det[\mathcal{S}w_{b-1}(s_a)]_{a,b=1,\ldots N}}{\det[\mathcal{S}w_{b-1}(a-1)]_{a,b=1,\ldots N}}\prod_{j=1}^N\frac{\mathcal{S}\omega(s_j)}{\mathcal{S}\omega(j-1)}
\end{equation}
for the first statement along the results~\eqref{spherical-polynomial} and~\eqref{spherical-Polya}. Pulling the factors of $\mathcal{S}\omega(s_j)$ and $\mathcal{S}\omega(j-1)$ into the respective determinants and employing the convolution formula~\eqref{U(1)_convolve}, we obtain the claim. Similarly, we can do it for the second claim of the proposition.

For the third claim we start from~\eqref{convolution-factorisation-fixed} and have for the spherical transform of $U=U_1VU_2V^\dagger\in\U(N)$ 
\begin{equation}
\mathcal{S}^{(U)}(s)=\left(\prod_{j=0}^{N-1}j!\right)\frac{\det[x_a^{s_b}]_{a,b=1\ldots,N}}{\Delta_N(x)\Delta_N(s)}\prod_{j=1}^N\frac{\mathcal{S}\omega(s_j)}{\mathcal{S}\omega(j-1)},
\end{equation}
cf., Eq.~\eqref{spherical-function}. Anew, we pull the factors $\mathcal{S}\omega(s_j)$ and $\mathcal{S}\omega(j-1)$ into the determinant and use $x_a^{s_b}\mathcal{S}\omega(s_b)=\mathcal{S}[\omega(z'/x_a)](s_b)$, this time. The bijectivity of the spherical transform concludes the proof.
\hfill$\square$

The last statement of Theorem~\ref{thm:jpdfs} can be also rewritten in terms of a Harish-Chandra-like group integral identity. 

\begin{corollary}[Group Integral Identity for P\'olya Ensembles]\label{cor:groupintegral}\
	
	Let $f_N^{(U)}$ be a unitarily invariant cyclic P\'olya ensemble on $\U(N)$ associated to the weight $\omega\in\widetilde{L}_N^1(\mathbb{S}_1)$. Then, it satisfies the group integral identity
	\begin{equation}\label{group-int-identity}
	\int_{\U(N)}d\mu(U) f_N^{(U)}(Uy^\dagger U^\dagger x)=\frac{1}{\prod_{j=0}^{N-1}\mathcal{S}\omega(j)}\frac{\det[\omega(x_a/y_b)]_{a,b=1,\ldots,N-1}}{\Delta_N(x^\dagger)\Delta_N(y)}
	\end{equation}
	for all non-degenerate $x,y\in\mathbb{S}_1^N$. For degenerate $x$ and/or $y$ one needs to apply l'H\^opital's rule.
\end{corollary}

Let us underline that this statement can be readily extended to non-positive functions instead of probability densities. The weight $\omega$ only needs to satisfy suitable integrability and differentiability.

{\bf Proof of Corollary~\ref{cor:groupintegral}:}

We can understand the integral~\eqref{group-int-identity} as a probability density in $x\in\mathbb{S}_1^N$ when multiplying it with the factor $|\Delta_N(x)|^2/N!$. Indeed, the function
\begin{equation}
\tilde{p}(x)=\frac{|\Delta_N(x)|^2}{N!}\int_{\U(N)}d\mu(U) f_N^{(U)}(Uy^\dagger U^\dagger x)
\end{equation}
is evidently non-negative and symmetric under permutations in the elements of $x=\diag(x_1,\ldots,x_N)$. It is normalised because of
\begin{equation}
\begin{split}
\int_{\mathbb{S}_1^N}\left(\prod_{j=1}^N\frac{dx_j}{2\pi i x_j}\right)\tilde{p}(x)=&\int_{\mathbb{S}_1^N}\left(\prod_{j=1}^N\frac{dx_j}{2\pi i x_j}\right)\frac{|\Delta_N(x)|^2}{N!}\int_{\U(N)}d\mu(U) f_N^{(U)}(y^\dagger U^\dagger xU)\\
\overset{V=U^\dagger xU}{=}&\int_{\U(N)}d\mu(V) f_N^{(U)}(y^\dagger V)\\
\overset{V\to yV}{=}&\int_{\U(N)}d\mu(V) f_N^{(U)}(V)=1.
\end{split}
\end{equation}
In the second equality we have used that the measure of the matrix $V=U^\dagger xU$ distributed along $\left(\prod_{j=1}^N dx_j/(2\pi i x_j)\right)|\Delta_N(x)|^2d\mu(U) /N!$ is again the normalised Haar measure $d\mu(V)$ on the unitary group $\U(N)$.

With this knowledge we can compute the spherical transform of $\tilde{p}(x)$ which is
\begin{equation}
\begin{split}
\mathcal{S}\tilde{p}(s)=&\int_{\mathbb{S}_1^N} \left(\prod_{j=1}^N\frac{dz_j}{2\pi i x_j}\right)\tilde{p}(x)\Phi(x;s)\\
=&\int_{\U(N)}d\mu(V) f_N^{(U)}(y^\dagger V)\Phi(V;s)\\
=&\int_{\U(N)}d\mu(V) f_N^{(U)}(V)\Phi(yV;s)\\
=&\int_{\U(N)}d\mu(V) \int_{\U(N)}d\mu(W)f_N^{(U)}(V)\Phi(yWVW^\dagger;s)\\
=&\mathcal{S}^{(U)}(s)\Phi(y;s).
\end{split}
\end{equation}
In the penultimate step, we have exploited the unitary invariance of the measure $ f_N^{(U)}(V)d\mu(V)$ and introduced a Haar distributed unitary matrix $W\in\U(N)$. The final line shows that $\mathcal{S}\tilde{p}(s)$ agrees with the spherical transform of the random matrix $VWyW^\dagger$ where $V$ is drawn from the distribution $f_N^{(U)}$. Comparison with~\eqref{fixed.jpdf} closes the proof.
\hfill$\square$

\begin{remark}[Laurent Series of the Weight]\
	
	Due to the $2\pi$ periodicity of the weight $\omega\in \widetilde{L}_N^1(\mathbb{S}_1)$, we can write it in terms of a Laurent series
	\begin{equation}\label{Laurent-series}
	\omega(z')=\sum_{s=-\infty}^\infty u_s {z'}^{-s}\ {\rm with}\ \mathcal{S}\omega(s)=u_s\in\mathbb{C}.
	\end{equation}
	The differentiability of $\omega$ on $\mathbb{S}_1$, has to be $(N-2)$-times continuously differentiable and $(N-1)$-times almost everywhere, and the integrability conditions have some consequences for the coefficients $|u_s|$. For instance, the condition $\partial^{N-1}\omega\in L^1(\mathbb{S}_1)$ implies that $|s^{N-1}u_s|$ is bounded for all $s\in\mathbb{Z}$, because
	\begin{equation}
		|s^{N-1}u_s|= \left|s^{N-1}\int_{\mathbb S_1} \omega(z)z^{s-1}\dv z\right|
		\le\int_{\mathbb S_1} \left|(z\partial_z)^{N-1}\omega(z)\frac{\dv z}{z}\right|
	\end{equation} 
	and the right hand side is bounded by a linear combination of the integral of $|\partial^{k}w|$ where $k=1,\ldots,N-1$. So $|u_s|$ is bounded from above by a constant times $|s|^{-N+1}$ for large $s$ so that the absolute convergence of the Laurent series is given at least on the complex unit circle for all $N\geq 3$ and does not require any regularisation such as a Gaussian in the limit of a diverging variance. Whether the Laurent series converges on a ring or is even entire depends on the explicit form of the Fourier coefficients $u_s$.
	
	Moreover, we would like to mention that the property $[\omega(z')]^*={z'}^{N-1}\omega(z')$ for $z'\in\mathbb{S}_1$ is equivalent to the relation
	\begin{equation}\label{coefficient-rel}
	u_s^*=u_{N-1-s}.
	\end{equation}
	Additionally, the positivity of the normalisation constant~\eqref{norm-polynomial} which is for cyclic polynomial ensembles equal to
	\begin{equation}\label{CN-Polya}
	C_N=\prod_{j=0}^{N-1}j!\mathcal{S}\omega(j)=\prod_{j=0}^{N-1}j!u_j>0
	\end{equation}
	implies that $u_s\neq0$ for all $s=0,\ldots, N-1$. Especially for odd $N=2M+1$, we even obtain that $u_M$ is a positive real number, as~\eqref{coefficient-rel} implies $C_N=M! u_M \prod_{j=0}^{M-1}j!(2M-j)!|u_j|^2$. Since the joint probability density is invariant under multiplying $\omega$ with a positive constant one can set $u_M=1$ for odd $N=2M+1$.
	
	For even $N=2M$, we have even the freedom to rescale the weight $\omega$ with a non-zero real number so that one can choose the coefficients $u_{M-1}=u_M^*$ to be a phase in a suitable complex half-plane. We will make use of that later in the proof of Theorem~\ref{thm:unique-Laurent},
	when showing that the Laurent series is unique up to a global normalisation factor for a P\'olya ensemble. This will be our next main result,
	already stated in the Introduction.
\end{remark}

{\bf Proof of Theorem~\ref{thm:unique-Laurent}:}

Let $p_N^{(1)}$ and $p_N^{(2)}$ be the joint probability densities that correspond to $\omega_1$ and $\omega_2$, respectively, and $\mathcal{S}\omega_1(s')=u_{s'}^{(1)}$ and $\mathcal{S}\omega_2(s')=u_{s'}^{(2)}$ be their Laurent coefficients. Our starting point had been that $p_N^{(1)}=p_N^{(2)}$ although the weights are different. The uniqueness, up to a normalisation constant, is based on the injectivity of the spherical transform which means
\begin{equation}\label{start-eq}
\prod_{j=1}^N\frac{u^{(1)}_{s_j}}{u^{(1)}_{j-1}}=\mathcal{S}p_N^{(1)}(s)=\mathcal{S}p_N^{(2)}(s)=\prod_{j=1}^N\frac{u^{(2)}_{s_j}}{u^{(2)}_{j-1}}
\end{equation}
for all $s=\diag(s_1,\ldots,s_N)\in\mathbb{I}_N$.

We choose an integer $s'\notin\{0,\ldots,N-1\}$, an $l\in\{0,\ldots,N-1\}$, $s_1=s'$ and $(s_2,\ldots,s_N)$ as a permutation of the set $\{0,\ldots,N-1\}\setminus\{l\}$. Then, almost all terms cancel in the ratios of~\eqref{start-eq} and it simplifies to
\begin{equation}\label{intermediate-eq}
\frac{u^{(1)}_{s'}}{u^{(1)}_{l}}=\frac{u^{(2)}_{s'}}{u^{(2)}_{l}}\ \Leftrightarrow\ u^{(1)}_{s'}=\frac{u^{(1)}_{l}}{u^{(2)}_{l}}u^{(2)}_{s'}.
\end{equation}
This equation holds for all integers $s'\notin\{0,\ldots,N-1\}$ and $l=0,\ldots,N-1$. 

\textit{\underline{For $N=2M+1$ odd.}} We take $l=M$ and define $C=u^{(1)}_{M}/u^{(2)}_{M}>0$, and all coefficients with $s'\notin\{0,\ldots,N-1\}$ are related in a unified way like $u^{(1)}_{s'}=Cu^{(2)}_{s'}$. In the last step, we choose $s'=s_0$, where we know that $u^{(1)}_{s_0}=Cu^{(2)}_{s_0}\neq0$, and $l\in\{0,\ldots,N-1\}$ anew arbitrary, which yields
\begin{equation}\label{intermediate-eq.b}
\frac{Cu^{(2)}_{s_0}}{u^{(1)}_{l}}=\frac{u^{(1)}_{s_0}}{u^{(1)}_{l}}=\frac{u^{(2)}_{s_0}}{u^{(2)}_{l}}\ \Leftrightarrow\ u^{(1)}_{l}=Cu^{(2)}_{l}.
\end{equation}
Combining this knowledge with the Laurent series representation of the weight we have $\omega_1(z')=C\omega_2(z')$ with $C=u^{(1)}_{M}/u^{(2)}_{M}>0$.

\textit{\underline{For $N=2M$ even.}} We choose $l=M-1$ and $l=M$ yielding the two equations
\begin{equation}
u^{(1)}_{s'}=\frac{u^{(1)}_{M}}{u^{(2)}_{M}}u^{(2)}_{s'}\ {\rm and}\ u^{(1)}_{s'}=\frac{u^{(1)}_{M-1}}{u^{(2)}_{M-1}}u^{(2)}_{s'}=\left(\frac{u^{(1)}_{M}}{u^{(2)}_{M}}\right)^*u^{(2)}_{s'}.
\end{equation}
Either $u^{(1)}_{s'}$ vanishes and so does $u^{(2)}_{s'}$ or we can divide both equations telling us that the phase
\begin{equation}
\frac{u^{(1)}_{M}}{u^{(2)}_{M}}\left(\frac{u^{(2)}_{M}}{u^{(1)}_{M}}\right)^*=1\ \Leftrightarrow\ \frac{u^{(1)}_{M}}{u^{(2)}_{M}}=\left(\frac{u^{(1)}_{M}}{u^{(2)}_{M}}\right)^*
\end{equation}
is unity. Defining $C=u^{(1)}_{M}/u^{(2)}_{M}\in\mathbb{R}\setminus\{0\}$, we obtain $u^{(1)}_{s'}=Cu^{(2)}_{s'}$ for any integer $s'\notin\{0,\ldots,N-1\}$. From here it works along the same lines as for odd $N$, which concludes the proof.
\hfill$\square$

Theorem~\ref{thm:unique-Laurent} is not as trivial as it looks. The condition of a non-vanishing Laurent coefficient $u_{s_0}\neq0$ for an integer $s_0\neq0,\ldots,N-1$ is crucial. Actually, it tells us that there is only one cyclic P\'olya ensemble which does not satisfy this condition and, hence, for which this proposition is not applicable, which is the Haar measure and will be discussed as the first example of a P\'olya ensemble.

\subsection{Examples for Cyclic P\'olya Ensembles}\label{sec:examples}

\subsubsection{The Haar Measure}\label{sec:Haar}

The Haar distributed unitary matrices build a P\'olya ensemble because of Corollary~\ref{cor:spher.Haar}. Equation~\eqref{spherical-Haar} can be used to backwards-engineer what the corresponding weight $\omega^{\rm (Haar)}$ is, i.e., we find the geometric sum
\begin{equation}\label{geom-Haar}
\omega(z')=\sum_{s=0}^{N-1}{z'}^{-s}=\frac{1-{z'}^{-N}}{1-{z'}^{-1}}.
\end{equation}
But as already pointed out before, this is not the only sum which leads to the Haar measure.

\begin{proposition}[Ambiguity of the Weight for the Haar Measure]\label{prop:cuep1}\
	
	Every weight of the form
	\begin{align}\label{Haar-sum}
	\omega(z')=\sum_{s=0}^{N-1}u_s{z'}^{-s}
	\end{align}
	with $u_{s}=u^\ast_{N-1-s}\neq0$, and $u_{(N-1)/2}>0$ if $N$ is odd, yields a cyclic P\'olya ensemble that is the Haar measure on $\U(N)$, in particular it gives the joint probability density function~\eqref{jpdf-Polya}.
\end{proposition}

{\bf Proof of Proposition~\ref{prop:cuep1}:}

Due to the bijectivity of the spherical transform we only need to show that Eq.~\eqref{spherical-Polya} is equal to~\eqref{spherical-Haar}. Certainly, because of $\mathcal{S}\omega(s')=u_{s'}$ with $u_{s'}$ being the Laurent coefficient we notice that the finite sum~\eqref{Haar-sum} yields that the indices $s_j$ are restricted to the interval $[0,N-1]$. Therefore, we have
\begin{equation}
\mathcal{S}\omega(s)=\prod_{j=1}^N\frac{u_{s_j}}{u_{j-1}}\chi_{[0,N-1]}(s_j)
\end{equation}
for all $s\in\mathbb{I}_N$. The set $\mathbb{I}_N$ implies pairwise different components in the multi-index $s$. However there are only $N$ integers in $[0,N-1]$ so that $(s_1,\ldots,s_N)$ has to be a permutation of $(0,\ldots,N-1)$. This guarantees for the product $\prod_{j=1}^N u_{s_j}/u_{j-1}=1$ and, thus, we are left with the product of the characteristic functions which is indeed Eq.~\eqref{spherical-Haar}, finishing the proof.
\hfill$\square$

One very suitable weight yielding the Haar measure which we will encounter later on is of a binomial form
\begin{equation}\label{binomial-Haar}
\omega_N^{\rm(Haar)}(z')=\sum_{j=0}^{N-1}\binom{N-1}{j}{z'}^{-j}=(1+{z'}^{-1})^{N-1}=2^{N-1}\left[\cos\left(\frac{\theta}{2}\right)\right]^{N-1}e^{-i(N-1)\theta/2}
\end{equation}
for $z=e^{i\theta}$ with $\theta\in]-\pi,\pi[$.

\begin{remark}
	
	(1)~(Stability of the Haar Measure) From the defining property of the Haar measure on $\U(N)$, we must have that for $U_1, U_2 \in\U(N)$ with
	$U_1$ Haar distributed, the product $U_1U_2\in\U(N)$ is also Haar distributed. This can be also be seen on the level of the spherical transform which is
	\begin{equation}
	\begin{split} 
	\mathcal{S}^{(U_1U_2)}(s)=\mathcal{S}^{(U_1)}(s)\mathcal{S}^{(U_2)}(s)=\mathcal{S}^{(U_2)}(s)\prod_{j=1}^N\chi_{[0,N-1]}(s_j)=\prod_{j=1}^N\chi_{[0,N-1]}(s_j)=\mathcal{S}^{(U_1)}(s).
	\end{split}
	\end{equation}
	In the second to last step, we have used that $(s_1,\ldots,s_N)$ has to be a permutation of $\{0,\ldots,N-1\}$ so that we can evaluate $\mathcal{S}^{(U_2)}(s)$ as $\mathcal{S}^{(U_2)}(0,\ldots,N-1)=1$ due to the normalisation. \\
	(2)~Haar distributed compact Lie groups can be considered as the oldest of the random matrix ensembles. The finding and parameterisation of the group invariant measures, in particular in relation to $\U(N)$, was a topic in mathematics~\cite{Hurwitz} when nobody thought about random matrix theory as an independent field; see the review~\cite{DF17}.
\end{remark}

\subsubsection{Brownian Motion on a Circle}\label{sec:Brownian}

In~\cite{Liechty}, the Dyson-Brownian motion on a circle has been considered, especially on the unitary group $\U(N)$. In particular, the heat equation
\begin{equation}\label{heat-equation}
\partial_t f_N(U;t)=\mathcal {L}_Uf_N(U,t)
\end{equation}
has been solved for some initial condition $f_N(U;0)$ and the Laplace-Beltrami operator $\mathcal L_U$ that corresponds to the unique (up to a normalisation) group invariant Haar metric on $\U(N)$ which also creates the Haar measure. If the initial condition is a Dirac delta function on $\U(N)$ at the point $U_0$, then, $f_N(U;t|U_0)$ describes the probability density of $U_t$.

The induced Laplace-Beltrami operator $\mathcal{L}_z$ for the eigenvalues $z=\diag(z_1,\ldots,z_N)\in\mathbb{S}_1^N$ of the matrix $U$ is explicitly given by
\begin{equation}
\mathcal L_z=\frac{1}{|\Delta_N(z)|}\left(\sum_{j=1}^Nz_j\partial_{z_j}z_j^*\partial_{z_j^*}\right)|\Delta_N(z)|=\frac{1}{\Delta_N(z^*)}\left(-\sum_{j=1}^N\left[\,z_j\partial_{z_j}+\frac{N-1}{2}\right]^2\right)\Delta_N(z^*).
\end{equation}
The fundamental solution $u(z;t)$ of the heat kernel is the initial boundary value problem
\begin{equation}\label{fund-kernel-heat}
\partial_t u(z;t)=\mathcal {L}_zu(z;t)\quad{\rm for}\ z\in\mathbb{S}_1^N\quad {\rm and}\quad u(z;0)=\prod_{j=1}^N\delta(z_j-1),
\end{equation}
where $\delta(z_j-1)$ is the Dirac delta function on the complex unit circle with the property 
\begin{equation}
\int_{\mathbb{S}_1}f(z')\delta(z'-z_0) dz'=f(z_0)
\end{equation}
for any $z_0\in\mathbb{S}_1$ and any function $f$ on $\mathbb{S}_1$. The Dirac delta functions in~\eqref{fund-kernel-heat} enforce that the initial point of the Brownian motion is at $U_0=\mathbf{1}_N$. The kernel $u(z;t)$ has been computed in~\cite{Liechty} and it is given in the following proposition.

\begin{proposition}[Proposition 1.1 in~\cite{Liechty}]\label{prop:heat}\
	
	The fundamental solution $u(z;t)$ of the heat equation~\eqref{fund-kernel-heat} times $|\Delta_N(z)|^2$ and a proper normalisation is a joint probability density of the eigenvalues of a cyclic P\'olya ensemble, which we call cyclic Gaussian ensemble, with the weight
	\begin{equation}\label{Gauss-Polya}
	\omega_N^{\rm (Gauss)}(z';t)=\sum_{s=-\infty}^\infty\exp\left[-t\left(s+\frac{1-N}{2}\right)^2\right]{z'}^{-s},
	\end{equation}
	which is a Jacobi-theta function~\cite[\S 20.2(i)]{NIST}.
	Especially, the joint probability density $p_N^{\rm(Gauss)}$ has the form~\eqref{jpdf-Polya}.
\end{proposition}

From the knowledge of the fundamental solution of the heat equation, we can deduce two simple consequences. By Corollary~\ref{cor:spherical-Polya} the cyclic Gaussian ensemble has the spherical transform
\begin{equation}
\s p_N^{\rm(Gauss)}(s;t)=\prod_{j=1}^N \exp\left[-t\left(s_j+\frac{1-N}{2}\right)^2+t\left(j-1+\frac{1-N}{2}\right)^2\right].
\end{equation}
We made use of this in our recent work~\cite{ZKF2019} and also introduced it in~\eqref{spherical-inverse} to regularise the inverse of the spherical transform.

The second consequence is yielded by Theorem~\ref{thm:jpdfs} part (3) implying the transition kernel of the heat equation when the initial condition is not $U_0=\mathbf{1}_N$ but an arbitrary $U_0$ with the eigenvalues $x=\diag(x_1,\ldots,x_N)\in\mathbb{S}_1^N$. Then, the distribution $p_N^{\rm(Gauss)}(y;t|x)$ of the eigenvalues $y=\diag(y_1,\ldots,_N)\in\mathbb{S}_1^N$ of $U_t$ is given as in~\eqref{fixed.jpdf} with $\omega$ being the Jacobi-theta function $\omega_N^{\rm(Gauss)}$, see~\eqref{Gauss-Polya}.

\subsubsection{The Circular Jacobi Ensemble}\label{sec:Jacobi}

As a third ensemble, we would like to mention the circular (or cyclic) Jacobi ensemble~\cite{WF2000,BD2002,FLT2020}, which has the joint probability density
\begin{equation}\label{cjepdf}
\begin{split}
p_N^{\rm(Jac)}(z;\alpha,\gamma)=&\frac{|\Delta_N(z)|^2}{\widetilde{C}_N}\prod_{j=1}^N\left|(1+z_j)^{\alpha-2i\gamma}\right|=2^{\alpha N}\frac{|\Delta_N(e^{i\theta})|^2}{\widetilde{C}_N }\prod_{j=1}^N\left[\cos\left(\frac{\theta_j}{2}\right)\right]^\alpha e^{\gamma\theta_j},
\end{split}
\end{equation}
where $\alpha>-1, \gamma\in\mathbb{R}$ are two parameters. 
This was introduced in (\ref{CJi}).
To render the square root taken on the right side meaningful, we assume that the cut is taken along the negative real half-axis meaning for the angles $\theta_j\in]-\pi,\pi[$ of the complex phases $z_j=e^{i\theta_j}\in\mathbb{S}_1$. Indeed, the point $z_j=-1$ is a Fisher-Hartwig singularity~\cite{ES1997} as the confining potential may even experience a jump of a finite height-difference when $\alpha=0$, meaning it can mimic a potential step. The asymptotic behaviour of the spectrum close to such a singular point is described by the confluent hypergeometric kernel, see~\cite{BD2002,FLT2020}.

The density~\eqref{cjepdf} has been considered in several works, for instance because of its relation to Selberg integrals~\cite[\S 3.9]{Forresterbook}. In~\cite{DIK2011,DIK2012,DIV2010}, the authors considered a broader class by choosing $\alpha,\,\gamma\in\mathbb{C}$ with ${\rm Re}(\alpha)>-1$. However, we would like to focus on probability weights.

To see that this ensemble is a cyclic P\'olya ensemble we rewrite the term
\begin{equation}
|1+z|^\alpha=z^{-\alpha/2}(1+z)^\alpha,\quad z\in\mathbb S_1
\end{equation}
and observe that
\begin{equation}
\begin{split}
\det\left[(-z_a\partial_a)^{b-1}z_a^{\nu}(1+z_a)^\mu\right]_{a,b=1,\ldots,N}=&\Delta_N(z^*)\prod_{j=1}^N\frac{\Gamma[\mu+1]}{\Gamma[\mu-j+2]}z_j^{\nu+N-1}(1+z_j)^{\mu-N+1}
\end{split}
\end{equation}
for any two exponents $\mu,\nu\in\mathbb{C}$. In this way, we can identify the weight (\ref{Jacobi-weight}) 
as that for the present P\'olya ensemble. When comparing this result with the weight~\eqref{binomial-Haar}, we recognise that the Haar measure is a very particular form of the cyclic Jacobi ensemble namely for $\alpha=\gamma=0$. Indeed, this could be expected from the joint probability density~\eqref{cjepdf}, so this is a good consistency check.

The spherical transform easily follows from
\begin{equation}\label{sphere.Jac}
\begin{split}
\mathcal{S}\omega_N^{\rm(Jac)}(s';\alpha,\gamma)=&\int_{-\pi}^\pi\frac{d\theta}{2\pi}e^{is'\theta}e^{-i\alpha\theta/2-i(N-1)\theta+\gamma\theta}(1+e^{i\theta})^{\alpha+N-1}\\
=&\frac{\Gamma[N+\alpha]}{\Gamma[N+\alpha/2-s'+i\gamma]\Gamma[\alpha/2+s'-i\gamma+1]},
\end{split}
\end{equation}
which is
\begin{equation}
\mathcal{S}p_N^{\rm(Jac)}(s;\alpha,\gamma)=\prod_{j=1}^N\frac{\Gamma[N+\alpha/2-j+i\gamma+1]\Gamma[\alpha/2+j-i\gamma]}{\Gamma[N+\alpha/2-s_j+i\gamma]\Gamma[\alpha/2+s_j-i\gamma+1]}.
\end{equation}
What has been elegantly carried out has been essentially a Selberg integral~\cite{Forresterbook}. This can be particularly seen for the Morris integral~\cite{Morris} which is the normalisation factor
\begin{equation}\label{Ho1}
\begin{split}
\widetilde{C}_N=&\int_{]-\pi,\pi[^N} |\Delta_n(e^{i\theta})|^2\prod_{j=1}^Ne^{\gamma \theta_j}|1+e^{i\theta_j}|^\alpha \frac{d\theta_j}{2\pi}\\
=&\left(\prod_{j=0}^{N-1}\frac{\Gamma[\alpha+N-j]}{\Gamma[\alpha+N]}\right)\int_{\mathbb{S}_1^N}\left(\prod_{j=1}^N\frac{dz_j}{2\pi i z_j}\right)\Delta_N(z^*)\det[(-z_a\partial_a)^{b-1}z_a^{-\alpha/2-i\gamma}(1+z_a)^{\alpha+N-1}]_{a,b=1,\ldots,N}\\
=&N!\left(\prod_{j=0}^{N-1}\frac{\Gamma[\alpha+N-j]}{\Gamma[\alpha+N]}\right)\left(\prod_{j=0}^{N-1}j!\mathcal{S}\omega_N^{\rm(Jac)}(j)\right)\\
=&\prod_{j=0}^{N-1}\frac{\Gamma(1+\alpha+j)\Gamma(j+2)}{|\Gamma(1+\alpha/2+i\gamma+j)|^2}.
\end{split}
\end{equation}

\subsubsection{Bilateral Hypergeometric Ensemble}\label{sec:Bilateral}

We have seen several examples of cyclic P\'olya ensemble. In fact, the weights~\eqref{geom-Haar} and~\eqref{Jacobi-weight} are very special cases of the \textit{bilateral hypergeometric series} (\ref{biletral-hyper}).
The full potential of the bilateral hypergeometric function unfolds when studying products of cyclic Jacobi ensembles or similar ensembles.
By forming the Fourier series with the help of (\ref{Ho1}) in the case $N = 1$, and
from the definition~\eqref{biletral-hyper}, we can identify the cyclic Jacobi weight~\eqref{Jacobi-weight} according to
\begin{equation}\label{1h1}
\begin{split}
\omega_N^{\rm(Jac)}(z';\alpha,\gamma)=&\sum_{s=-\infty}^\infty\frac{\Gamma[N+\alpha]}{\Gamma[N+\alpha/2-s+i\gamma]\Gamma[\alpha/2+s-i\gamma+1]} {z'}^{-s}\\
=&\frac{\Gamma[N+\alpha/2+i\gamma]\Gamma[N+\alpha]}{\Gamma[\alpha/2-i\gamma+1]} {_1}H_1\left[\left.
\begin{array}{c}
-N-\alpha/2-i\gamma+1\\
\alpha/2-i\gamma+1
\end{array}\right|-z'\right],
\end{split}
\end{equation}
where we have exploited Euler's reflection formula for the ratio $\Gamma[-N-\alpha/2-i\gamma+1]/\Gamma[N+\alpha/2-s+i\gamma]$. Additionally, we know from Theorem~\ref{thm:jpdfs}.2 that the product of two or more circular Jacobi matrices is still a P\'olya ensemble with a weight function which is equal to the convolution of all the weight functions. For instance, for the product $U=U_1U_2\cdots U_k$ where $U_j\in\U(N)$ is drawn from a cyclic Jacobi ensemble with the weight $\omega_j(z')=\omega^{\rm(Jac)}(z';\alpha_j,\gamma_j)$ the new P\'olya ensemble is associated to the weight
\begin{equation}
\begin{split}
\omega(z')=&\omega_1\ast\omega_2\ast\cdots\ast\omega_k(z')=\left(\prod_{l=1}^k\frac{\Gamma[N+\alpha_l/2+i\gamma_l]\Gamma[N+\alpha_l]}{\Gamma[\alpha_l/2-i\gamma_l+1]}\right){_k}H_k\left[\left.
\begin{array}{c}
a_1,\ldots,a_k\\
b_1,\ldots,b_k
\end{array}\right|-z'\right]
\end{split}
\end{equation}
with $a_j=-N-\alpha_j/2-i\gamma_j+1$ and $b_j=\alpha_j/2-i\gamma_j+1$. The weight stays a bilateral hypergeometric function only with more indices
in keeping with the general relation (\ref{HH}). As remarked in the paragraph containing (\ref{HH}), this shows that cyclic P\'olya ensembles with bilateral hypergeometric weights play the role of P\'olya ensembles with Meijer-G function weights in the study of unitary invariant products of positive definite Hermitian matrices.
Corollary~\ref{cor:spherical-Polya} part (2) shows that also the inverse random matrix $U^{-1}$ of bilateral hypergeometric random matrix remains in the class. Indeed in~\cite[Eq.~(6.1.1.4)]{Slater}, we can read off the identity
\begin{align}\label{reciprocal}
\phq\left[\left.
\begin{aligned}
a_1,~\ldots,~a_p\\
b_1,~\ldots,~b_q
\end{aligned}\right|z'
\right]={_q}H_p\left[\left.
\begin{aligned}
1-b_1,~\ldots,~1-b_q\\
1-a_1,~\ldots,~1-a_p
\end{aligned}\right|{z'}^{-1}\right],
\end{align}
which readily shows this claim.

\subsubsection{Constructing Cyclic P\'olya Ensembles from Rank $1$ Multiplications}\label{sec:Rank1}

Many P\'olya ensembles can be created in a very simple way via multiplying specific exponentiated rank-$1$ unitary matrices. This will be shown in the ensuing paragraphs.

\begin{definition}[Cyclic Rank-1 Jacobi Ensemble]\label{def:rank1}\
	
	Let $\gamma\in\mathbb{R}$. A cyclic rank-1 Jacobi matrix $U_\gamma\in\U(N)$ is a random matrix which can be decomposed like
	\begin{equation}\label{rank1-Jacobi}
	U_\gamma=V\diag(\mathbf{1}_{N-1},-x)V^\dagger
	\end{equation}
	with a complex phase $x=e^{i\theta}\in\mathbb{S}_1$ distributed by the density
	\begin{equation}
	p_\gamma(x)=\frac{|\Gamma[(N+1)/2+i\gamma]|^2}{(N-1)!}\left|(1+x)^{N-1-2i\gamma}\right|=2^{N-1}\frac{|\Gamma[(N+1)/2+i\gamma]|^2}{(N-1)!}\left[\cos\left(\frac{\theta}{2}\right)\right]^{N-1}e^{\gamma \theta}
	\end{equation}
	with $\theta\in]-\pi,\pi[$ and $V\in\U(N)$ a Haar distributed unitary matrix. We denote the set of these matrices by $\mathcal{R}_1(N)$.
\end{definition}

The chosen name of these ensembles becomes clear when comparing it with the joint probability density~\eqref{cjepdf}.

Its spherical transform is the first we will compute as it is the starting point of constructing cyclic P\'olya ensembles.

\begin{proposition}[Spherical Transform of Cyclic Rank-1 Jacobi Matrices]\label{prop:rank1-spher}\
	
	Let $\gamma\in\mathbb{R}$. The spherical transform of a random matrix $U_\gamma\in\mathcal{R}_1(N)$ is
	\begin{equation}\label{rank1-spher}
	\mathcal{S}^{(U_\gamma)}(s)=\prod_{j=1}^{N}\frac{(1-N)/2-i\gamma+j-1}{(1-N)/2-i\gamma+s_j}
	\end{equation}
	for all $s=\diag(s_1,\ldots,s_N)\in\mathbb{I}_N$.
	The case for $\gamma=0$ and odd $N$ has to be understood via l'H\^opital's rule.
\end{proposition}

{\bf Proof of Proposition~\ref{prop:rank1-spher}:}

Choosing a $U_\gamma=V\diag(1,\ldots,1,-x)V^\dagger\in\mathcal{R}_1(N)$ with a $\gamma\neq0$, we perform the integral
\begin{equation}
\begin{split}
\mathcal{S}^{(U_\gamma)}=&\frac{|\Gamma[(N+1)/2+i\gamma]|^2}{(N-1)!}\int_{\mathbb{S}_1}\frac{dx}{2\pi i x}\left|(1+x)^{N-1-2i\gamma}\right|\int_{\U(N)} d\mu(V) \Phi(V\diag(\mathbf{1}_{N-1},-x)V^\dagger;s)\\
=&\frac{|\Gamma[(N+1)/2+i\gamma]|^2}{(N-1)!}\int_{\mathbb{S}_1}\frac{dx}{2\pi i x}(1+x)^{N-1}x^{(1-N)/2-i\gamma}\Phi(\diag(\mathbf{1}_{N-1},-x);s).
\end{split}
\end{equation}
In the last step, we have exploited the unitary invariance of the spherical function $\Phi$. We apply l'H\^opital's rule in~\eqref{spherical-function} to find
\begin{equation}
\Phi(\diag(\mathbf{1}_{N-1},-x);s)=(N-1)!\frac{\det[ \begin{array}{c|c} s_a^{b-1} & (-x)^{s_a}\end{array}]_{\substack{a=1\ldots,N \\ b=1,\ldots,N-1}}}{(1+x)^{N-1}\Delta_N(s)}.
\end{equation}
The factor $(1+x)^{N-1}$ cancels, and the integral over $x$ can be carried out by using
\begin{equation}
\int_{\mathbb{S}_1}\frac{dx}{2\pi i x}x^{(1-N)/2-i\gamma+s_a}=\int_{-\pi}^\pi\frac{d\theta}{2\pi}e^{i\theta((1-N)/2-i\gamma+s_a)}=(-1)^{s_a}\frac{\sin(\pi[(1-N)/2-i\gamma])}{\pi[(1-N)/2-i\gamma+s_a]}.
\end{equation}
Thus, the sign $(-1)^{s_a}$ is cancelling and we are left with
\begin{equation}
\det\left[ \begin{array}{c|c} s_a^{b-1} & \displaystyle\frac{1}{(1-N)/2-i\gamma+s_a} \end{array}\right]_{\substack{a=1\ldots,N \\ b=1,\ldots,N-1}}=(-1)^{N-1}\frac{\Delta_N(s)}{\prod_{j=1}^N((1-N)/2-i\gamma+s_j)}.
\end{equation}
This kind of determinant has been employed in several other works such as in~\cite{KG2010,BF1994}. Collecting everything, we arrive at the claim~\eqref{rank1-spher}. The case $\gamma=0$ can be found via the limit $\gamma\to0$ which works out too for odd $N$ as then the numerator and denominator in~\eqref{rank1-spher} vanish like $\gamma$.
\hfill$\square$

From Eq.~\eqref{rank1-spher}, we see that the cyclic rank-1 Jacobi ensembles are essentially P\'olya ensembles if we do not care that the joint probability density has to be a function but can be a general distribution. The distribution shows itself in the $N-1$ fixed eigenvalues of $U_\gamma$ at $1$. The corresponding weight is
\begin{equation}\label{weight-rank1}
\begin{split}
\omega_N^{\rm(rank)}(z';\gamma)=&\lim_{t\to0^+}\sum_{s=-\infty}^\infty\frac{-i\tilde{\gamma}}{(1-N)/2-i\gamma+s}{z'}^{-s}\exp\left[-ts^2\right]\\
=&\frac{\tilde{\gamma}}{2\sinh(\pi[\gamma+i(1-N)/2])}(-z')^{(1-N)/2-i\gamma}
\end{split}
\end{equation}
with $\tilde{\gamma}=\gamma$ when $N$ is odd and $\tilde{\gamma}=1$ when it is even. This sum can be computed with the help of Poisson's summation rule. The combination $(-z')$ ensures the that the cut and, hence, the jump of the weight is along the positive real axis of $z'$ as the root has the cut commonly along the negative one. It is the reason why this weight is not differentiable and how it creates the $N-1$ eigenvalues at $z'=1$ when one interprets the weight as a distribution. The non-analyticity at $z'=1$ also guarantees the $2\pi$-periodicity of the weight.

The form of the weight~\eqref{weight-rank1}
upon comparison with~\eqref{Jacobi-weight}
 is the reason for the name rank-1 Jacobi ensembles. P\'olya ensembles with no degenerate eigenvalues can be created by multiplying at least $N$ cyclic rank-1 Jacobi matrices.

\begin{corollary}[Cyclic P\'olya Ensembles from Cyclic Rank-1 Jacobi Matrices]\label{cor:cpe-rank1jac}\
	
	Let $L\geq N$ be a positive integer, $\gamma_1,\ldots,\gamma_L\in\mathbb{R}$ be real constants, and $U_{\gamma_1},\ldots,U_{\gamma_L}\in\mathcal{R}_1(N)$ be cyclic rank-1 Jacobi matrices. Then, the product matrix $U=U_{\gamma_1}U_{\gamma_2}\cdots U_{\gamma_L}$ is equivalent in distribution with a random matrix drawn from a cyclic P\'olya ensemble associated to the weight
	\begin{equation}
	\omega(z')={_L}H_{L}\left[\left.
	\begin{array}{c}
	\frac{1-N}{2}-i\gamma_1,\ldots,\frac{1-N}{2}-i\gamma_L\\
	\frac{3-N}{2}-i\gamma_1+1,\ldots,\frac{3-N}{2}-i\gamma_L+1
	\end{array}\right|z\right].
	\end{equation}
\end{corollary}

{\bf Proof of Corollary~\ref{cor:cpe-rank1jac}:}

Due to the factorisation of the spherical transform $\mathcal{S}^{(U)}=\prod_{j=1}^N\mathcal{S}^{(U_{\gamma_j})}$, see~\eqref{convolution-factorisation-random}, we can identify the coefficients of the bilateral hypergeometric function because of the relation $1/a=\Gamma[a]/\Gamma[a+1]$ for any $a\neq0$. The differentiability and integrability of $\omega$, see~\eqref{L1-set.b}, follows from the absolute convergence of the series as the modules of the coefficients drop off like $1/|s|^L$ for $|s|\to\infty$. The non-negativity of the joint probability density follows from the fact that $U$ is a product of random matrices and that the convolution of probability measures stay probability measures. This closes the proof.
\hfill$\square$

\begin{remark}[Generation of Gamma Functions in the Spherical Transform]\
	
	When taking infinite products, we can even generate Gamma functions in the Laurent series via the Weierstrass formula
	\begin{equation}
	\Gamma[x+1]=e^{-\gamma_{\rm E}x}\prod_{l=1}^\infty\frac{\exp[x/l]}{1+x/l}
	\end{equation}
	with $\gamma_{\rm E}\approx 0.58$ the Euler-Mascheroni constant. For instance, when defining the unitary matrices $V_l=e^{i/l}U_{l+\nu}$ with $U_{l+\nu}\in\mathcal{R}_1(N)$ for $l\in\mathbb{N}$ and $\nu>-1$, their spherical transform is equal to
	\begin{equation}
	\mathcal{S}^{(V_l)}(s)=\prod_{j=1}^{N}e^{i(s_j-j+1)/l}\frac{(1-N)/2-i(l+\nu)+j-1}{(1-N)/2-i(l+\nu)+s_j}.
	\end{equation}
	Thus, we find for the infinite product $V=e^{-i\gamma_{\rm E}}V_1V_2\cdots$ the spherical transform
	\begin{equation}\label{spher.Gamma}
	\mathcal{S}^{(V)}(s)=\lim_{L\to\infty}e^{-i\gamma_{\rm E}\sum_{j=1}^N(s_j-j+1)}\prod_{l=1}^L\mathcal{S}^{(V_l)}(s)=\prod_{j=1}^{N}\frac{\Gamma[\nu+1+i(s+[1-N]/2)]}{\Gamma[\nu+1+i(j-[1+N]/2)]}.
	\end{equation}
	The corresponding cyclic P\'olya ensemble yielding this spherical transform is the counterpart of the Laguerre (induced Ginibre) ensemble~\cite{AI2015,KK2016,KK2019} for the multiplicative convolution on ${\rm GL}_{\mathbb{C}}(N)$ and a Muttalib-Borodin ensemble~\cite{FKK2017}, where the weight function is the Gumble distribution times an exponential factor $e^{-\nu x}$, for the additive convolution on the Hermitian matrices. Thence, we coin the corresponding weight as
	\begin{equation}\label{weight.Gin}
	\omega_N^{\rm (Gin)}(z';\nu)=\sum_{s=-\infty}^\infty\Gamma[\nu+1+i(s+[1-N]/2)]{z'}^{-s}
	\end{equation}
	and call the corresponding ensemble the cyclic Ginibre ensemble.
	The limit~\eqref{spher.Gamma} can be indeed carried over to the probability density level as the corresponding series of the inverse transform~\eqref{spherical-inverse-jpdf} is absolutely convergent when $l\geq 2N$.
	
	As a side remark, we have exploited the fact that the multiplication of a unitary random matrix $U\in\U(N)$ with a constant phase $z_0\in\mathbb{S}_1$ results in the spherical transform
	\begin{equation}\label{spher.const.phase}
	\mathcal{S}^{(z_0U)}(s)=\mathcal{S}^{(U)}(s)\prod_{j=1}^Nz_0^{s_j-j+1}.
	\end{equation}
	This can be readily checked by the definitions~\eqref{spherical-transform} and~\eqref{spherical-function}. 
\end{remark}

\subsection{The Positivity and the Relation to Cyclic P\'olya Frequency Functions}\label{sec:CPF}

As we have learned, we can write the weights of cyclic P\'olya ensembles in terms of Laurent series. The problem is that not any Laurent series satisfies the requirement that the probability density of the eigenvalues is non-negative. To solve this hurdle we consider P\'olya frequency functions on $\mathbb{S}_1$.

\begin{definition}[P\'olya Frequency Functions on $\mathbb{S}_1$]\label{def:Polyafunction}\
	\begin{enumerate}
		\item	Let $N=2M+1\in2\mathbb{N}+1$ be odd. Then, a function $g:\mathbb{S}_1\mapsto\mathbb{R}_+$ satisfying
		\begin{equation}\label{Pol-frequ-odd}
		\frac{\Delta_{2m+1}(x)\Delta_{2m+1}(y^{-1})}{[\det (xy^{-1})]^{m}}\det\left[g(x_ay_b^{-1})\right]_{a,b=1,\ldots,2m+1}\geq0,
		\end{equation}
		for all $x,y\in\mathbb{S}_1^{2m+1}$ and $m=0,1,\ldots,M$, is called P\'olya frequency function of order $2M+1$ (see~\cite{Karlin,KRS1994}). Here $\det\left[xy^{-1}\right]:=\prod_{j=1}^Nx_jy_j^{-1}$
		\item	Let $N=2M\in2\mathbb{N}$ be even. Then, a function $g:\mathbb{S}_1\to\mathbb{C}$ satisfying $[g(z)]^*=zg(z)$ and
		\begin{equation}\label{Pol-frequ-even}
		\frac{\Delta_{2m}(x)\Delta_{2m}(y^{-1})}{[\det (xy^{-1})]^{m-1}}\det\left[g(x_ay_b^{-1})\right]_{a,b=1,\ldots,2m}\geq0,
		\end{equation}
		for all $x,y\in\mathbb{S}_1^{2m}$ and $m=1,\ldots,M$, is called P\'olya frequency function of order $2M$.
	\end{enumerate}
\end{definition}

P\'olya frequency functions for odd orders $N=2M+1$ have been already defined in~\cite{Karlin,KRS1994}, while in~\cite[Ch 9]{Karlin} the above definition is instead referred to as the extended cyclic P\'olya frequency function of order $2M+1$. The subtle difference of the definition for odd and even dimensions is born out the complex conjugation of the Vandermonde determinant, see~\eqref{Vand-ident}. This is also the reason why the function $g$ needs to be complex. Certainly, the condition $[g(z)]^*=zg(z)$ only means that $(z)^{1/2}g(z)$ is real if we cut the complex plane along the negative real axis.

\begin{example}\label{ex:Polya.frequ}\
	
	Let us give some examples of such cyclic P\'olya frequency functions.
	\begin{enumerate}
		\item	The function
		\begin{equation}\label{Haar-Polya}			
		g_N^{\rm(Haar)}(z')=\left\{\begin{array}{cl} \displaystyle \sum_{j=0}^{2M}\binom{2M}{j}(z')^{M-j}=2^{2M}\left[\cos\left(\frac{\theta}{2}\right)\right]^{2M},\ & N=2M+1,\\ \displaystyle \sum_{j=0}^{2M-1}\binom{2M-1}{j}(z')^{M-1-j}=2^{2M-1}\left[\cos\left(\frac{\theta}{2}\right)\right]^{2M-1}e^{-i\theta/2},\ & N=2M, \end{array}\right.
		\end{equation}
		with $z'=e^{i\theta}\in\mathbb{S}_1$ with $\theta\in[-\pi,\pi]$ is a cyclic P\'olya frequency function of order $N$, respectively whether $N$ is odd or even. Note that we need to cut the complex plane along the negative real axis to match the two ends when $N=2M$ is even. The $N=2M+1$ case is referred to as the \textbf{De la Valle\'e Poussin kernel} in~\cite[Ch 9 \S 3]{Karlin}, and a proof that such kernel is indeed a cyclic P\'olya frequency function of order $2M+1$ can be also seen in~\cite[Ch 9 Thm 3.1]{Karlin}.
		
		The superscript is reminiscent of the weight for the Haar measure. Indeed, we have
		\begin{equation}
			 \omega_{2M+1}^{\rm(Haar)}(z')={z'}^{-M}g_{2M+1}^{\rm(Haar)}(z'),\quad \omega_{2M}^{\rm(Haar)}(z')={z'}^{1-M}g_{2M}^{\rm(Haar)}(z').
		\end{equation}
		
		The property of the cyclic P\'olya frequency function follows from the the group integral~\eqref{group-int-identity} and noticing that
		\begin{equation}
		\omega_{2M+1}^{\rm(Jac)}(z';2M-2m,0)={z'}^{m}g_{2M+1}^{\rm(Haar)}(z'),\quad \omega_{2M}^{\rm(Jac)}(z';2M-2m,0))={z'}^{m-1}g_{2M}^{\rm(Haar)}(z')	
		\end{equation}  
		are the weights of cyclic Jacobi ensembles for any $m\leq M$ which is known to create a random matrix ensemble and thus its probability density is positive on the left hand side of Eq.~\eqref{group-int-identity}.
		
		Along the same lines one can show that the functions related to the general cyclic Jacobi weights,
		\begin{equation}\label{Jac-Polya}			
		g_N^{\rm(Jac)}(z';\alpha,\gamma)=\left\{\begin{array}{cl} \displaystyle \sum_{j=-\infty}^{\infty}\frac{(z')^{-j}}{\Gamma[M+\alpha/2-j+i\gamma+1]\Gamma[M+\alpha/2+j-i\gamma+1]},\ & N=2M+1,\\ \displaystyle \sum_{j=-\infty}^{\infty}\frac{(z')^{-j}}{\Gamma[M+\alpha/2-j+i\gamma+1]\Gamma[M+\alpha/2+j-i\gamma]},\ & N=2M, \end{array}\right.
		\end{equation}
		are cyclic P\'olya frequency functions of order $N+\lceil\alpha\rceil$, where $\lceil.\rceil$ is the ceil function yielding the smallest integer which is larger than or equal to $\alpha$.
		
		\item	Let $N=2M+\chi$ with $\chi=0,1$, encoding whether $N$ is even or odd. With the help of the group integral~\eqref{group-int-identity}, one can also show that the Jacobi-theta function
		\begin{equation}\label{Gauss-Polya-frequ}			
		g_{2-\chi}^{\rm(Gauss)}(z';t)=\left\{\begin{array}{cl} \displaystyle \sum_{j=-\infty}^{\infty}\exp\left[-tj^2\right](z')^{-j},\ & \chi=1,\\ \displaystyle \sum_{j=-\infty}^{\infty}\exp\left[-t\left(j-\frac{1}{2}\right)^2\right](z')^{-j},\ & \chi=0, \end{array}\right.
		\end{equation}
		is a cyclic P\'olya frequency function. This time it is of infinite odd or even order, respectively, as we can create the cyclic Gaussian weight for any dimension via 
		\begin{equation}
			\omega_{2M+1}^{\rm(Gauss)}(z';t)={z'}^{M} g_{1}^{\rm(Gauss)}(z';t),\quad \omega_{2M}^{\rm(Gauss)}(z';t)={z'}^{M-1} g_{2}^{\rm(Gauss)}(z';t).
		\end{equation}
		
		\item	Also the weight for the rank-1 Jacobi matrices can be related with the following cyclic P\'olya frequency functions
		\begin{equation}\label{Rank-Polya}
		g_{2-\chi}^{\rm(rank)}(z';\gamma)=\left\{\begin{array}{cl} \displaystyle (-z')^{-i\gamma}=e^{\gamma\theta}[\cosh(\gamma\pi)-\sinh(\gamma\pi){\rm sign}(\theta)],\ & \chi=1,\\ \displaystyle -i(-z')^{-1/2-i\gamma}=e^{(\gamma-i/2)\theta}[\sinh(\gamma\pi)+\cosh(\gamma\pi){\rm sign}(\theta)],\ & \chi=0, \end{array}\right.
		\end{equation}
		where $z'=e^{i\theta}\in\mathbb{S}_1$ with $\theta\in]-2\pi,2\pi[$ and ${\rm sign}(\alpha)$ yields the sign of $\alpha\in\mathbb{R}\setminus\{0\}$ and vanishes when $\alpha=0$, with the relation
		\begin{equation}
			\omega_{2M+1}^{\rm(rank)}(z';\gamma)={z'}^{-M} g_{1}^{\rm(rank)}(z';\gamma),\quad 
			\omega_{2M}^{\rm(rank)}(z';\gamma)={z'}^{1-M} g_{2}^{\rm(rank)}(z';\gamma).
		\end{equation} 
		We have chosen the interval $]-2\pi,2\pi[$ instead of $[0,2\pi[$ to prove that the above function is indeed a cyclic P\'olya frequency function, the reason being we encounter differences $\theta_a-\phi_b$ of two angles $\theta_a,\phi_b\in[0,2\pi[$ when choosing the phases $x_a=e^{i\theta_a}$ and $y_b=e^{i\phi_b}$ in~\eqref{Pol-frequ-odd} and~\eqref{Pol-frequ-even}. We prove this in the following proposition.
	\end{enumerate}
\end{example}

\begin{proposition}[Cyclic P\'olya Frequency Function of Rank-1 Case]\label{prop:Rank-1.Polya}\
	
	The functions~\eqref{Rank-Polya} are P\'olya frequency functions of any odd or even order, respectively.
\end{proposition}

{\bf Proof of Proposition~\ref{prop:Rank-1.Polya}:}

To check the statement for the odd dimensional case, we compute
\begin{equation}\label{rank-1-check.odd}
\begin{split}
&\frac{\Delta_{2m+1}(x)\Delta_{2m+1}(y^{-1})}{[\det (xy^{-1})]^{m}}\det\left[g_1^{\rm(rank)}(x_ay_b^{-1};\gamma)\right]_{a,b=1,\ldots,2m+1}\\
=&\left(\prod_{1\leq k<l<2m+1}4\sin\left[\frac{\theta_l-\theta_k}{2}\right]\sin\left[\frac{\phi_l-\phi_k}{2}\right]\right)\left(\prod_{l=1}^{2m+1}e^{\gamma(\theta_l-\phi_l)}\cosh(\gamma\pi)\right)\\
&\times\det[1-\tanh(\gamma\pi){\rm sign}(\theta_a-\phi_b)]_{a,b=1,\ldots,2m+1}
\end{split}
\end{equation}
for $\gamma\neq0$, as it trivially vanishes for $\gamma=0$.
We do not loose generality when assuming an ordering of the angles as follows $0\leq \theta_1<\theta_2<\ldots< \theta_{2m+1}<2\pi$ and $0\leq \phi_1<\phi_2<\ldots< \phi_{2m+1}<2\pi$. Indeed, the determinant and the sine functions are zero whenever $\theta_l=\theta_k$ or $\phi_l=\phi_k$ for some $l\neq k$. Moreover, their product is symmetric under permutation of the angles $\{\theta_j\}_{j=1,\ldots,2m+1}$ as well as of the angles $\{\phi_j\}_{j=1,\ldots,2m+1}$. Additionally, one can show that whenever two angles $\phi_l\leq\theta_k<\theta_{k+1}\leq\phi_{l+1}$ or $\theta_l\leq\phi_k<\phi_{k+1}\leq\theta_{l+1}$ for some $l,k\in\{1,\ldots,2m+1\}$ with $\theta_{2m+2}=\theta_{1}+2\pi$ and $\phi_{2m+2}=\phi_{1}+2\pi$ the remaining determinant vanishes as either two rows or two columns become exactly the same. Therefore, the two sets of angles can only satisfy one of the two possible interlacing conditions
\begin{equation}\label{interlacing}
0\leq\theta_1\leq\phi_1\leq\theta_2\leq\ldots\leq\theta_{2m+1}\leq\phi_{2m+1}<2\pi\ {\rm or}\ 0\leq\phi_1\leq\theta_1\leq\phi_2\leq\ldots\leq\phi_{2m+1}\leq\theta_{2m+1}<2\pi.
\end{equation}
Due to the symmetry in the two sets of angles we can assume $0\leq\theta_1\leq\phi_1\leq\ldots$ 
This implies that the matrix $T\in\mathbb{R}^{(2m+1)\times(2m+1)}$ with the entries $T_{ab}=1-\tanh(\gamma\pi){\rm sign}(\theta_a-\phi_b)$ is explicitly $T_{ab}=1-\tanh(\gamma\pi)$ and $T_{ba}=1+\tanh(\gamma\pi)$ for all $a>b$ and on the diagonal we have $T_{aa}\in\{1,1+\tanh(\gamma\pi)\}$. Whenever there is a $k\in\{1,\ldots,2m+1\}$ with $T_{kk}=1+\tanh(\gamma\pi)$ we can subtract the last $2m-k+2$ columns with the $k$th and the first $k-1$ columns with $[1-\tanh(\gamma\pi)]/[1+\tanh(\gamma\pi)]$ times the $k$-th one. Thus the determinant of $T$ evaluates to $\det(T)=[1+\tanh(\gamma\pi)][\tanh(\gamma\pi)]^{2m}\prod_{j\neq k}[1-{\rm sign}(\theta_j-\phi_j)]\geq0$. If all diagonal entries are $T_{aa}=1$, the determinant becomes $\det[T]=[\tanh(\gamma\pi)]^{2m}\geq0$ as can be readily checked by induction in the dimension $m$. Plugging this insight into~\eqref{rank-1-check.odd} shows our claim for the odd dimensional case.

Similarly, we approach the even dimensional case where we have
\begin{equation}\label{rank-1-check.even}
\begin{split}
&\frac{\Delta_{2m}(x)\Delta_{2m}(y^{-1})}{[\det (xy^{-1})]^{m-1}}\det\left[g_2^{\rm(rank)}(x_ay_b^{-1};\gamma)\right]_{a,b=1,\ldots,2m}\\
=&\left(\prod_{1\leq k<l<2m}4\sin\left[\frac{\theta_l-\theta_k}{2}\right]\sin\left[\frac{\phi_l-\phi_k}{2}\right]\right)\left(\prod_{l=1}^{2m}e^{\gamma(\theta_l-\phi_l)}\cosh(\gamma\pi)\right)\\
&\times\det\left[\tanh(\gamma\pi)+{\rm sign}(\theta_a-\phi_b)\right]_{a,b=1,\ldots,2m}.
\end{split}
\end{equation}
From this expression we can anew read off that we can order the angles without loss of generality and that the interlacing condition is again valid, i.e.,
\begin{equation}\label{interlacing.b}
0\leq\theta_1\leq\phi_1\leq\theta_2\leq\ldots\leq\theta_{2m}\leq\phi_{2m}<2\pi\ {\rm or}\ 0\leq\phi_1\leq\theta_1\leq\phi_2\leq\ldots\leq\phi_{2m}\leq\theta_{2m}<2\pi.
\end{equation}
Since the symmetry between the two sets of angles allows us to choose $0\leq\theta_1\leq\phi_1\leq\ldots$, we consider the determinant of the matrix $T\in\mathbb{R}^{2m\times2m}$ with the entries $T_{ab}=\tanh(\gamma\pi)+{\rm sign}(\theta_a-\phi_b)$. This time this means $T_{ab}=\tanh(\gamma\pi)+1$ and $T_{ba}=\tanh(\gamma\pi)-1$ for all $a>b$ and on the diagonal we have $T_{aa}\in\{\tanh(\gamma\pi),\tanh(\gamma\pi)-1\}$. Whenever there is a $k\in\{1,\ldots,2m\}$ with $T_{kk}=\tanh(\gamma\pi)-1$, we subtract the last $2m-k+1$ columns with the $k$th one and the first $k-1$ columns with the $k$-th one times $[\tanh(\gamma\pi)+1]/[\tanh(\gamma\pi)-1]$ which leads to 
\begin{equation}
	\det[T]=[1-\tanh(\gamma\pi)]\prod_{j\neq k}[1-{\rm sign}(\theta_j-\phi_j)]\geq0.	
\end{equation}
For the case of all diagonal entries $T_{aa}=\tanh(\gamma\pi)$, we find $\det[T]=1>0$. Plugging this into~\eqref{interlacing.b} closes the proof for the even dimensional case.
\hfill$\square$

With the Examples~\ref{ex:Polya.frequ}, we can obtain already a big class of cyclic P\'olya frequency functions, namely via the multiplicative convolution on the complex sphere. This is the analogue of the convolution theorem for P\'olya frequency functions on the real line, see~\cite[Prop 7.1.5]{Karlin}, and for the odd case it is also implicitly implied in~\cite[Ch. 9 Thm. 4.1]{Karlin})

\begin{proposition}[Convolution of Cyclic P\'olya Frequency Functions]\label{prop:conv.Polya}\
	
	Let $g_1$ and $g_2$ be two cyclic P\'olya frequency functions of order $N$ and suitably integrable so that $g_1(\widetilde{z})g_2(z'\widetilde{z}^{-1})$ is absolutely integrable in $\widetilde{z}\in\mathbb{S}_1$ for all $z'\in\mathbb{S}_1$ with respect to the Haar measure on $\mathbb{S}_1$. Then, the convolution
	\begin{equation}\label{conv.Polya}
	g_1\ast g_2(z')=\int_{\mathbb{S}_1}\frac{d\widetilde{z}}{2\pi i \widetilde{z}}g_1(\widetilde{z})g_2(z'\widetilde{z}^{-1})
	\end{equation}
	is also a cyclic P\'olya frequency functions of order $N$.
\end{proposition}

{\bf Proof of Proposition~\ref{prop:conv.Polya}:}

Again let $N=2M+\chi$ with $\chi=0,1$. Then, the reality condition can be readily checked
\begin{equation}
\begin{split}
[g_1\ast g_2(z')]^*=&\int_{\mathbb{S}_1}\frac{d\widetilde{z}}{2\pi i \widetilde{z}}[g_1(\widetilde{z})]^*[g_2(z'\widetilde{z}^{-1})]^*\\
=&\int_{\mathbb{S}_1}\frac{d\widetilde{z}}{2\pi i \widetilde{z}}\widetilde{z}^{1-\chi} g_1(\widetilde{z}) \left(\frac{z'}{\widetilde{z}}\right)^{1-\chi} g_2(z'\widetilde{z}^{-1})={z'}^{1-\chi} g_1\ast g_2(z').
\end{split}
\end{equation}

Next, we choose two sets of phases $x,y\in\mathbb{S}_1^{2m+\chi}$ and $m=1,\ldots,M$ and compute
\begin{equation}
\begin{split}
&\frac{\Delta_{2m+\chi}(x)\Delta_{2m+\chi}(y^{-1})}{[\det (xy^{-1})]^{m+\chi-1}}\det\left[g_1\ast g_2(x_ay_b^{-1})\right]_{a,b=1,\ldots,2m+\chi}\\
=&\frac{\Delta_{2m+\chi}(x)\Delta_{2m+\chi}(y^{-1})}{[\det (xy^{-1})]^{m+\chi-1}}\det\left[\int_{\mathbb{S}_1}\frac{d\widetilde{z}}{2\pi i \widetilde{z}}g_1(\widetilde{z}y_b^{-1})g_2(x_a\widetilde{z}^{-1})\right]_{a,b=1,\ldots,2m+\chi}\\
=&\frac{1}{(2m+\chi)!}\int_{\mathbb{S}_1^{2m+\chi}}\left(\prod_{j=1}^{2m+\chi}\frac{dz_j}{2\pi i z_j}\right)\frac{\Delta_{2m+\chi}(x)\Delta_{2m+\chi}(z^{-1})\det\left[g_2(x_az_b^{-1})\right]_{a,b=1,\ldots,2m+\chi}}{[\det (xz^{-1})]^{m+\chi-1}}\\
&\times\frac{\Delta_{2m+\chi}(z)\Delta_{2m+\chi}(y^{-1})\det\left[g_1(z_ay_b^{-1})\right]_{a,b=1,\ldots,2m+\chi}}{[\det (zy^{-1})]^{m+\chi-1}}\frac{1}{|\Delta_{2m+\chi}(z)|^2}\geq0.
\end{split}
\end{equation}
In the penultimate step, we have employed the Andr\'eief identity~\eqref{Andreief} with $k=0$ and $N=2m+\chi$, and the last inequality follows from the fact that the two functions are cyclic P\'olya frequency functions. This closes the proof.
\hfill$\square$

The two Propositions~\ref{prop:Rank-1.Polya} and~\ref{prop:conv.Polya} give rise to many other P\'olya frequency functions in a very constructive way. For instance the cyclic Ginibre case leads to the functions
\begin{equation}\label{Gin-Polya}			
g_{2-\chi}^{\rm(Gin)}(z';\nu)=\left\{\begin{array}{cl} \displaystyle \sum_{j=-\infty}^{\infty}\Gamma[\nu+1+ij](z')^{-j},\ & \chi=1,\\ \displaystyle \sum_{j=-\infty}^{\infty}\Gamma[\nu+1+i(j-1/2)](z')^{-j},\ & \chi=0, \end{array}\right.
\end{equation}
which are cyclic P\'olya ensembles of infinite odd or even order. However, the classification of cyclic P\'olya frequency functions is still incomplete with these example. It is reasonable to conjecture that the product of cyclic rank-1 Jacobi functions~\eqref{Rank-Polya}, the cyclic Gaussian function~\eqref{Gauss-Polya-frequ} and the multiplication with a constant phase $z_0$ may yield all cyclic P\'olya frequency functions of infinite order, as it is the case for their counterpart on the real line~\cite{Schoenberg}. A proof of this claim is still an open problem~\cite[p 384]{Schoenberg2}. Even worse is the situation of the classification of the cyclic P\'olya frequency functions at finite order. For instance, the cyclic Jacobi ensemble and here in particular the Haar measure of the unitary group $\U(N)$ are the only P\'olya frequency functions of a finite order which have been widely studied. One can expect that there are many more which fall into this class.

Undoubtedly there is a deep relation between cyclic P\'olya frequency function and cyclic P\'olya ensemble, which is the reason why we have named the ensembles in this way. It is evident, for instance, that a suitably differentiable and integrable cyclic P\'olya frequency function gives rise to a respective ensemble on $\U(N)$. The inverse statement is not so trivial as one needs to check that when $\omega_N(z')$ gives rise to a cyclic P\'olya ensemble on $\U(N)$ that $\omega_{N-2m}(z')={z'}^{m}\omega_N(z')$ corresponds to one on $\U(N-2m)$ for any $2m<N$. Although our examples presented above corroborate this statement, a
general proof is yet to be found. We content ourselves instead with the following theorem, which will be the last one in the present section.

\begin{theorem}[Relation of Cyclic P\'olya Ensembles and Frequency Functions]\label{thm:Polya.rel}\
	
	Let $g_\chi\in\widetilde{L}^1(\mathbb{S}_1)$ be $(2M+\chi-1)$-times differentiable and a cyclic P\'olya frequency function of odd ($\chi=1$) or even ($\chi=0$) order. Then, $\omega_{2M+\chi}(z')={z'}^{-M-\chi+1}g_\chi(z')$ is a weight associated to a cyclic P\'olya ensemble on $\U(2M+\chi)$.
\end{theorem}

{\bf Proof of Theorem~\ref{thm:Polya.rel}:}

The integrability and differentiability conditions stay the same when multiplying $g_\chi$ with the analytic phase factor ${z'}^{-M-\chi+1}$. The identity $[\omega_{2M+\chi}(z')]^*=z^{2M+\chi-1}\omega_{2M+\chi}(z')$ results from the pre-factor and the realness condition of the cyclic P\'olya frequency function $[g_\chi(z')]^*={z'}^{1-\chi}g_\chi(z')$. Thus, we need to prove the positivity of the joint probability density. For this aim, we divide either Eq.~\eqref{Pol-frequ-odd} or~\eqref{Pol-frequ-even} by $|\Delta_{2M+\chi}(y)|^2$ and take the limit from a non-degenerate $y$, say $y_a=\exp[i\epsilon a]$ with $\epsilon\to0$, to $y=\mathbf{1}_{2M+\chi}$ via l'H\^opital's rule. What we obtain is the joint probability density~\eqref{jpdf-Polya} up a normalisation constant. This density is indeed positive as it has been the case for any non-degenerate $y$ and $a=1,\ldots,2M+\chi$. This shows our claim.
\hfill$\square$

\subsection{Relationship with the Derivative Principle}

As already commented in the Introduction, there are various recent studies on P\'olya ensembles in many other matrix spaces including Hemitian matrix space and positive definite Hermtian matrix space~\cite{KK2016,KK2019,FKK2017,KR2016}. One usually introduces those classes of ensembles by giving exact formulae for the eigenvalue distribution, similarly to the present $\U(N)$ case~\eqref{jpdf-Polya}. All these representations have similar forms (see~\cite{FKK2017}) in terms of a product of Vandermonde determinant and another determinant with derivatives acting on a weight function $w$. The viewpoint taken is that such a structure with two determinants gives a determinantal point process and, hence, allows a study using a bi-orthogonal system to explicitly write down its correlation kernel~\cite{KK2019,Kieburg2019}.

From another viewpoint, Ref.~\cite{ZK2020} shows that matrices in those spaces with a certain group invariance have eigenvalue distributions with a similar structure, assuming only modest analytical requirements. In particular, for $U\in\U(N)$ being invariant under unitary conjugation, there exists a symmetric function $g:\mathbb S^{N}\mapsto \R$ such that
\begin{equation}\label{DP}
	p_N^{(U)}(z)=\frac{1}{\prod_{j=1}^Nj!}\Delta(z)\prod_{a<b}\left(z_a\partial_{a}-z_b\partial_{b}\right)g(z_1,\ldots,z_N)
\end{equation}
This is referred to as the \textbf{derivative principle}. The weight function $g$ is also unique, under some modest analytical requirements (which we believe can be relaxed by using distribution theory), as well as the requirement
\begin{equation}\label{DP_uniqueness}
	\int_{\mathbb S^N}g(z_1,\ldots,z_N)\prod_{j=0}^N\frac{z_j^{s_j}\,\dv z_j}{2\pi i z_j}\,=0\ {\rm with}\ s_1,\ldots,s_N\in\mathbb{Z},
\end{equation}
whenever $s_j=s_k$ for some $j\neq k\in\{1,2,\ldots,N\}$. It is also given as an existence theorem which shows a way to construct such a weight function as an average on $\U(N)$, by using a parametrisation of the unitary group~\cite[Appendix B]{ZK2020}.

Comparing~\eqref{DP} and~\eqref{jpdf-Polya}, we immediately notice that a P\'olya ensemble is obtained when the weight function $g$ is replaced by a product of univariate weight functions $w$ up to a scalar (the functions in the product must be identical because of the symmetry of $f$). This is however not possible as~\eqref{DP_uniqueness} cannot be met for any non-zero weight. Yet, one can add any homogeneous solution $g_H$ of the differential equation $\prod_{a<b}\left(z_a\partial_{a}-z_b\partial_{b}\right)g_H(z_1,\ldots,z_N)=0$ to $g$ so that $g+g_H$ is such a product of $w$. Therefore we say that the structure of P\'olya ensemble is a natural choice for determinantal point processes on $\U(N)$.

Let us also compare the two expressions of the Haar measure in~\cite{ZK2020} and Proposition~\ref{prop:cuep1}. It is given in~\cite{ZK2020} that
\begin{equation}\label{pdfHaar1}
	p^{(\text{Haar})}_N(z)=\frac{1}{\prod_{j=0}^Nj!}\Delta(z)\prod_{a<b}\left(z_a\partial_{a}-z_b\partial_{b}\right)\sum_{\rho\in S_N}\prod_{j=1}^Nz_j^{-(\rho(j)-1)},
\end{equation}
while Proposition~\ref{prop:cuep1} gives
\begin{equation}\label{pdfHaar2}
	p^{(\text{Haar})}_N(z)=\frac{1}{\prod_{j=0}^Nj!}\Delta(z)\prod_{a<b}\left(z_a\partial_{a}-z_b\partial_{b}\right)\frac{1}{\prod_{j=0}^{N-1}u_j}\prod_{j=1}^N\sum_{k=0}^{N-1}u_sz_j^{-s},
\end{equation}
with $u_s$ satisfying the conditions given below~\eqref{Haar-sum}. It can be checked that these two expressions are equivalent. Thus notice that any function of the form $h(z_jz_k)$ for any $j,k=0,\ldots,N-1$ is a homogeneous solution of the differential equation $\prod_{a<b}\left(z_a\partial_{a}-z_b\partial_{b}\right)g_H(z_1,\ldots,z_N)=0$. So after expanding the product in~\eqref{pdfHaar2}, the only monomials surviving the action of the Vandermonde differential operator are such that no two $z_j$ and $z_k$ would have the same power. As the highest power of a $z_j$ is $N-1$, only the monomials $z_1^{-(\rho(1)-1)}z_2^{-(\rho(2)-1)}\ldots z_{N}^{-(\rho(N)-1)}$ for some permutation $\rho\in S_N$ are surviving the derivative operator $\prod_{a<b}\left(z_a\partial_{a}-z_b\partial_{b}\right)$. Summing over those permutations gives exactly~\eqref{pdfHaar1}.

\section{Eigenvalue Statistics of Products of Unitary Random Matrices}\label{sec:prod}

In this section we derive the kernels of cyclic P\'olya ensembles (subsection~\ref{sec:EVstat.Polya}) and products of these ensembles with either fixed matrices (subsection~\ref{sec:EVstat.fixed}) or cyclic polynomial ensemble (subsection~\ref{sec:EVstat.polynomial}). We especially aim at simple formula in terms of bi-orthonormal functions. Here, we adapt the approach and notions of~\cite{Kieburg2019}.

\begin{definition}[Bi-orthonormal Pair of Functions]\label{def:biorthonormal}\
	
	A set $\{(P_{j},Q_j)\}_{j=0,\ldots,N-1}$ is said to be a bi-orthonormal pair of functions of a cyclic polynomial ensemble associated to the weights $\{w_{j}\}_{j=0,\ldots,N-1}\subset L_N^1(\mathbb{S}_1)$ if the following three properties are satisfied:
	\begin{enumerate}
		\item	the linear span of polynomials is ${\rm span}_{j=0,\ldots,N-1}\{P_j\}={\rm span}_{j=0,\ldots,N-1}\{z^j\}$,
		\item	the linear span of weights is ${\rm span}_{j=0,\ldots,N-1}\{Q_j\}={\rm span}_{j=0,\ldots,N-1}\{w_j\}$,
		\item	for any $a,b=0,\ldots,N-1$ we have $\int_{\mathbb{S}_1} [dz'/(2\pi i z')]P_a(z')Q_b(z')=\delta_{ab}$.
	\end{enumerate}
\end{definition}

With the aid of a bi-orthonormal pair of functions the kernel of a determinantal point process, cf., Eq.~\eqref{k-point}, takes a very compact form, namely~\cite{Borodin}
\begin{equation}\label{kernel-general}
K_N(z_1,z_2)=\sum_{j=0}^{N-1}P_j(z_1)Q_j(z_2).
\end{equation}
One reason why we are interested in constructing such functions, and for obtaining corresponding structured forms of the correlation kernel, is for future use to compute asymptotics required in the analysis of double scaling limits. The utility of such developments is well evidenced for other matrix convolutions, e.g., see~\cite{AI2015,Forresterbook} and references therein. 

One last remark is in order. Evidently, a bi-orthonormal pair of functions is not uniquely given for a specific polynomial ensemble. One could fix this ambiguity by choosing $P_j$ to be a monic polynomial of order $j$. 
However this comes at the price of cumbersome normalisation constants, so is not adopted below.

\subsection{Eigenvalue Statistics of a Cyclic P\'olya Ensemble}\label{sec:EVstat.Polya}

As our first ensemble, we consider a cyclic P\'olya ensemble. A helpful quantity for the
 computation of the orthonormal pair of functions
 is the set $\mathbb{J}_l=\{0,1,\ldots,l-1\}$ for $l>0$ and $\mathbb{J}_0=\emptyset$ the empty set, as well as its complement $\mathbb{J}_l^{\rm c}=\mathbb{Z}\setminus\mathbb{J}_l$. We need too the ratio of gamma functions formula
\begin{equation}
\frac{\Gamma[N-j]}{\Gamma[-j]}=(-1)^{N-1}\frac{\Gamma[j+1]}{\Gamma[j-N+1]}
\end{equation} 
if $j$ is an integer which is larger than or equal to $N$. This allows us to write the result in a compact form.

\begin{proposition}[Kernel of a Cyclic P\'olya Ensemble]\label{prop:cyc.Pol}\
	
	A bi-orthonormal pair of functions $\{(P_{j},Q_j)\}_{j=0,\ldots,N-1}$ of the cyclic polynomial ensemble associated to the weight $\omega\in\widetilde{L}_N^1(\mathbb{S}_1)$ is
	\begin{equation}\label{biorth-cyc.Pol}
	\begin{split}
	P_j(z_1)=&\sum_{k\in\mathbb{J}_{j+1}}\frac{1}{(j-k)!k!}\frac{(-z_1)^k}{\mathcal{S}\omega(k)},\\
	Q_j(z_2)=&z_2\partial_2^j{z_2}^{j-1}\omega(z_2)=\lim_{t\to0}\sum_{l\in\mathbb{J}_{j}^{\rm c}}\frac{\Gamma[j-l]}{\Gamma[-l]}\mathcal{S}\omega(l)z_2^{-l}e^{-t (l+1-N)l}
	\end{split}
	\end{equation}
	for $j=0,\ldots,N-1$. The kernel is then the double sum
	\begin{equation}\label{kernel-cyc.Pol}
	\begin{split}
	K_N(z_1,z_2)=&\sum_{k\in\mathbb{J}_N}(z_1z_2^{-1})^k+\lim_{t\to0}\sum_{k\in\mathbb{J}_N}\sum_{l\in\mathbb{J}_N^{\rm c}}\frac{\Gamma[N-l]}{\Gamma[-l]\Gamma[N-k]\Gamma[k+1]}\frac{\mathcal{S}\omega(l)}{\mathcal{S}\omega(k)}\frac{(-z_1)^{k}z_2^{-l}}{k-l}e^{-t (l+1-N)l}.
	\end{split}
	\end{equation}
\end{proposition}

We underline that the formulas for the polynomials and weights imply very simple recurrence relations,
\begin{equation}
(j-z_1\partial_1)P_j(z_1)=P_{j-1}(z_1)\quad{\rm and}\quad (j+z_2\partial_2)Q_j(z_2)=Q_{j+1}(z_2).
\end{equation}
Thus, the differential operators in front of the bi-orthonormal functions can be understood as ladder operators and the formula of $Q_j(z_2)$ in terms of a differential operator is essentially a Rodrigues formula.

{\bf Proof of Proposition~\ref{prop:cyc.Pol}:}

The functions $Q_j$ are in the span of the weights $\{(-z'\partial)^j\omega(z')\}_{j=0,\ldots,N-1}$ because of the identity
\begin{equation}
z'\partial^j{z'}^{j-1}\omega(z')=\prod_{l=0}^{j-1}(z'\partial+l)\omega(z').
\end{equation}
Moreover, they are linearly independent which can be seen when computing their Fourier transform on $\mathbb{S}_1$ and using that $\mathcal{S}\omega(s)$ is at $N$ different points, namely at $s=0,\ldots,N-1$, non-vanishing.
The second identity in~\eqref{kernel-cyc.Pol} follows from the Laurent series representation of the weight.

The bi-orthonormality can be readily checked via direct computation. For this aim, we perform an integration by parts which is allowed as $\omega$ is $(N-2)$-times continuous differentiable and $2\pi$-periodic. Thus, we find
\begin{equation}
\begin{split}
\int_{\mathbb{S}_1} \frac{dz'}{2\pi i z'}P_a(z')Q_b(z')=&(-1)^b\int_{\mathbb{S}_1} \frac{dz'}{2\pi i z'}\sum_{k\in\mathbb{J}_{a+1}}\frac{(-1)^b}{\Gamma[k-b+1](a-k)!}\frac{(-z')^k}{\mathcal{S}\omega(k)}\omega(z')\\
=&\sum_{k\in\mathbb{J}_{a+1}}\frac{(-1)^{k-b}}{\Gamma[k-b+1](a-k)!}=\delta_{ab}.
\end{split}
\end{equation}
Here, we have used that $1/\Gamma[x+1]$ has zeros at negative integers so that all summands for $k<b$ are vanishing. This implies that the sum is zero whenever $a<b$. For $a>b$, we obtain a binomial sum yielding $(a-b)!(1-1)^{a-b}=0$, and for $a=b$ the sum only consists of the term $k=a=b$ rendering it equal to $1$.

For the kernel~\eqref{kernel-cyc.Pol}, we start from
\begin{equation}
\begin{split}
K_N(z_1,z_2)=&\sum_{j=0}^{N-1}P_j(z_1)Q_j(z_2)\\
=&\lim_{t\to0}\sum_{j=0}^{N-1}\sum_{k=0}^{j}\sum_{l=-\infty}^{-1}\frac{(j-l-1)!}{(j-k)!k!(-l-1)!}\frac{\mathcal{S}\omega(l)}{\mathcal{S}\omega(k)}(-z_1)^kz_2^{-l}e^{-t (l+1-N)l}\\
&+\lim_{t\to0}\sum_{j=0}^{N-1}\sum_{k=0}^{j}\sum_{l=j}^{\infty}(-1)^j\frac{l!}{(j-k)!k!(l-j)!}\frac{\mathcal{S}\omega(l)}{\mathcal{S}\omega(k)}(-z_1)^kz_2^{-l}e^{-t (l+1-N)l}\\
=&\lim_{t\to0}\sum_{k=0}^{N-1}\sum_{l=-\infty}^{-1}\left(\sum_{j=k}^{N-1}\frac{(j-l-1)!}{(j-k)!}\right)\frac{1}{k!(-l-1)!}\frac{\mathcal{S}\omega(l)}{\mathcal{S}\omega(k)}(-z_1)^kz_2^{-l}e^{-t (l+1-N)l}\\
&+\lim_{t\to0}\sum_{k=0}^{j}\sum_{l=0}^{\infty}\left(\sum_{j=k}^{\min\{N-1,l\}}\frac{(-1)^j}{(j-k)!(l-j)!}\right)\frac{l!}{k!}\frac{\mathcal{S}\omega(l)}{\mathcal{S}\omega(k)}(-z_1)^kz_2^{-l}e^{-t (l+1-N)l}.
\end{split}
\end{equation}
The Gaussian regularisation allows us to interchange the sums as they are all absolutely convergent. The sum over $j$ can be done via telescopic sums of the form
\begin{equation}
\begin{split}
(k-l)\sum_{j=k}^{N-1}\frac{(j-l-1)!}{(j-k)!}=&\frac{(N-l-1)!}{(N-k-1)!}=\frac{\Gamma(N-l)}{\Gamma(N-k)},\quad{\rm for}\ l<0,\\
(k-l)\sum_{j=k}^{N-1}\frac{(-1)^j}{(j-k)!(l-j)!}=&\frac{(-1)^{N}}{(N-k-1)!(l-N)!}=\frac{\Gamma(N-l)}{l!\Gamma(-l)\Gamma(N-k)},\quad{\rm for}\ l\geq N,
\end{split}
\end{equation}
and by the binomial sum for $l=0,\ldots,N-1$
\begin{equation}
\sum_{j=k}^{l}\frac{(-1)^j}{(j-k)!(l-j)!}=(-1)^k\delta_{lk}.
\end{equation}
Note that the latter sum is by definition zero when $l<k$. Putting everything together we find~\eqref{kernel-cyc.Pol}.
\hfill$\square$

The kernel can be cast into a simpler form of a one-fold integral as it has been done in sums and products with the following formula
\begin{equation}\label{integral-1ok}
\int_{0}^{2\pi}\frac{d\varphi}{2\pi}\, i\varphi\, e^{i k\varphi}=\frac{1}{k}\quad{\rm for}\ k\in\mathbb{Z}\setminus\{0\}.
\end{equation}
This yields a Christoffel-Darboux-like formula.

\begin{corollary}[Christoffel-Darboux-like Formula]\label{cor:Christ-Polya}\
	
	The kernel of the P\'olya ensemble of Proposition~\ref{prop:cyc.Pol} can be rewritten into the form
	\begin{equation}\label{kernel-cyc.Pol.b}
	K_N(z_1,z_2)=P_{N-1}(z_1)Q_{N-1}(z_2)+i\,\int_{0}^{2\pi}\frac{d\varphi}{2\pi}\varphi\, P_{N-2}(z_1e^{i \varphi})Q_{N-1}(z_2e^{i \varphi})+\frac{1-(z_1z_2^{-1})^{N-1}}{1-z_1z_2^{-1}}.
	\end{equation}
	When the weight satisfies $\omega\in \widetilde{L}_{N+1}^1(\mathbb{S}_1)$ this can be further reduced to (\ref{kernel-cyc.Pol.c}). Note that in the latter
	$Q_{N}(z')=(z'\partial+N-1)Q_{N-1}(z')$.
	\end{corollary}

Note, that for the Haar measure in~\eqref{kernel-cyc.Pol.c} we have $Q_N(z_2)=0$ because we take then the $N$-th derivative of a polynomial of order $N-1$.

{\bf Proof of Corollary~\ref{cor:Christ-Polya}:}

We only need to check that 
\begin{equation}
K_{N-1}(z_1,z_2)=i\,\int_{0}^{2\pi}\frac{d\varphi}{2\pi}\varphi\, P_{N-2}(z_1e^{i \varphi})Q_{N-1}(z_2e^{i \varphi})+\frac{1-(z_1z_2^{-1})^{N-1}}{1-z_1z_2^{-1}}
\end{equation}
as Eq.~\eqref{kernel-cyc.Pol.b} follows from $K_N(z_1,z_2)=P_{N-1}(z_1)Q_{N-1}(z_2)+K_{N-1}(z_1,z_2)$ and Eq.~\eqref{kernel-cyc.Pol.c} from the step $N-1\to N$. Essentially, we need only to argue that the integral~\eqref{integral-1ok} for $k\to k-l$ in~\eqref{kernel-cyc.Pol} can be interchanged with the sums. We underline that the regularisation can be omitted for $K_{N-1}(z_1,z_2)$ as then the summands drop off at least like $1/l^2$ because of the $(N-2)$-times continuous differentiability of the weight $\omega$. The interchange with the sum, then, results from the absolute integrability and convergence of the series leading to the desired form.
\hfill$\square$

Let us illustrate the results with the help of the cyclic Jacobi ensemble from subsection~\ref{sec:Jacobi}. The polynomials and weights are in this case equal to 
\begin{equation}\label{biorth-cyc.Jac}
\begin{split}
P_j(z_1;\alpha,\gamma)=&\sum_{k=0}^{\infty}\frac{\Gamma[N+\alpha/2-k+i\gamma]\Gamma[\alpha/2+k-i\gamma+1]}{\Gamma[j-k+1]\Gamma[N+\alpha]}\frac{(-z_1)^k}{k!}\\
=&\frac{\Gamma[N+\alpha/2+i\gamma]\Gamma[\alpha/2-i\gamma+1]}{j!\Gamma[N+\alpha]}{_2}F_1\left[\left.
\begin{array}{c}
-j\,,\,1+\alpha/2-i\gamma\\
1-N-\alpha/2-i\gamma
\end{array}\right|-z_1
\right]\\
=&\frac{(N+\alpha)\Gamma[\alpha/2-i\gamma+1]\Gamma[N-j+\alpha/2+i\gamma]}{j!\Gamma[N-j+\alpha+1]}{_2}F_1\left[\left.
\begin{array}{c}
-j\,,\,1+\alpha/2-i\gamma\\
N-j+\alpha+1
\end{array}\right|1+z_1
\right],\\
Q_j(z_2;\alpha,\gamma)=&z_2\partial_2^jz_2^{j-\alpha/2-i\gamma-N}(1+z_2)^{\alpha+N-1}\\
=&\sum_{l=0}^j\binom{j}{l}\frac{\Gamma[j-\alpha/2-i\gamma-N+1]\Gamma[\alpha+N]}{\Gamma[j-\alpha/2-i\gamma-N+1-l]\Gamma[\alpha+N-j+l]}z_2^{j-\alpha/2-i\gamma-N+1-l}(1+z_2)^{\alpha+N-1-j+l}\\
=&\frac{\Gamma[N+\alpha]}{\Gamma[N-j+\alpha]}\ |(1+z_2)^{\alpha-2i\gamma}|(1+z_2^*)^{N-1-j}{_2}F_1\left[\left.
\begin{array}{c}
-j\,,\,N-j+\alpha/2+i\gamma\\
N-j+\alpha
\end{array}\right|1+z_2^*
\right],
\end{split}
\end{equation}
cf., Eq.~\eqref{sphere.Jac} with $\alpha>-1$ and $\gamma\in\mathbb{R}$. For the polynomials, we have employed~\cite[Eq. (15.8.7)]{NIST}, and note that the infinite sum expression of $P_j$ is actually truncated to $k=j$, because the coefficients of the remaining terms vanishes. Hence, both sets of functions are expressible in terms of the hypergeometric functions~\cite[Eq. (16.2.1)]{NIST}
\begin{align}
{_p}F_q\left(
\left.\begin{array}{c}
a_1,\cdots,a_p\\
b_1,\cdots,b_q\end{array}\right| x
\right)=\frac{\prod_{j=1}^q\Gamma[b_j]}{\prod_{j=1}^p\Gamma[a_j]}\sum_{l=0}^\infty\frac{\prod_{j=1}^p\Gamma[a_j+l]}{\prod_{j=1}^q\Gamma[b_j+l]}\frac{x^l}{l!}.
\end{align}
The polynomial $P_j$ is comparable to the Routh-Romanovski polynomial appearing in the same ensemble in the work~\cite{FLT2020}.

In the case of the Haar-measure ($\alpha=\gamma=0$, see~\eqref{binomial-Haar}), we obtain highly non-trivial bi-orthonormal functions instead of the usually employed monomials. In the light of this, one may ask why we go through a more complicated expression. Here, we would like to emphasise that the results above hold for all cyclic P\'olya ensembles on $\U(N)$ and not only for the Haar measure. This has not been possible before without the technique outlined by us.

\subsection{Eigenvalue Statistics of a Product comprising a Fixed Matrix}\label{sec:EVstat.fixed}

Next we want to study the eigenvalue statistics of a product $U=U_1U_2$ of a cyclic P\'olya random matrix $U_2\in\U(N)$ that is associated to a weight $\omega\in\widetilde{L}_N^1(\mathbb{S}_1)$ and with a fixed unitary matrix $U_1\in\U(N)$. As we have seen in Theorem~\ref{thm:jpdfs} part (3) and in the proof of Corollary~\ref{cor:groupintegral}, the eigenvalue statistics of $U$ is not affected by whether the eigenvectors of $U_1$ are also fixed or randomly distributed as $U_2$ is unitarily invariant. What matters are only the eigenvalues $x=\diag(x_1,\ldots,x_N)\in\mathbb{S}_1^N$ of $U_1$.

Before we come to the bi-orthonormal pair of functions corresponding to the polynomial ensembles that is given by $U=U_1U_2$, we need to introduce the polynomial
\begin{equation}\label{chi-poly}
\chi_\omega(z')=\sum_{l=0}^{N-1}\frac{{z'}^l}{\mathcal{S}\omega(l)}.
\end{equation}
A similar polynomial has already been exploited in~\cite{Kieburg2019}.

\begin{proposition}[Kernel of a Cyclic P\'olya Ensemble times a Fixed Matrix]\label{prop:cyc.Pol.fixed}\
	
	Considering the setting of Theorem~\ref{thm:jpdfs} part (3), especially that the eigenvalues $x=\diag(x_1,\ldots,x_N)\in\mathbb{S}_1^N$ of $U_1$ are pairwise-different, the bi-orthonormal pair of functions $\{(P_{j},Q_j)\}_{j=0,\ldots,N-1}$ that describe the eigenvalue statistics of $U=U_1U_2$ are given by
	\begin{equation}\label{biorth-cyc.Pol.fixed}
	\begin{split}
	P_j(z_1)=&\int_{\mathbb{S}_1}\frac{dz'}{2\pi z'}\chi_\omega({z'}^{-1})\prod_{\substack{l=1,\ldots,N\\l\neq j+1}}\frac{z'z_1-x_l}{x_{j+1}-x_l},\quad Q_j(z_2)=\omega\left(\frac{z_2}{x_{j+1}}\right)
	\end{split}
	\end{equation}
	for $j=0,\ldots,N-1$. Assuming that the Laurent series of $\omega$ converges in a ring containing the complex unit circle, the kernel simplifies to the form of a double contour integral
	\begin{equation}\label{kernel-cyc.Pol.fixed}
	\begin{split}
	K_N(z_1,z_2)=&\int_{\mathbb{S}_1}\frac{dz'_1}{2\pi z'_1}\int_{\mathcal{C}}\frac{dz'_2}{2\pi z'_2}\frac{\chi_\omega(R^{-1}{z'}_1^{-1})\omega\left(z_2{z'}_2^{-1}\right)}{Rz'_1-z'_2}\prod_{l=1}^N\frac{R{z'}_1z_1-x_l}{z'_2-x_l},
	\end{split}
	\end{equation}
	where we choose a radius $R>1$ and a contour $\mathcal{C}$ encircling all eigenvalues $x=\diag(x_1,\ldots,x_N)\in\mathbb{S}_1^N$ counter-clockwise and close enough such that $|z_2|<R$ and stays in the ring of convergence of the Laurent series of $\omega$.
\end{proposition}

{\bf Proof of Proposition~\ref{prop:cyc.Pol.fixed}:}

The bi-orthonormality follows from the double contour integral identity
\begin{equation}\label{poly-int-ident}
\int_{\mathbb{S}_1}\frac{d\widetilde{z}}{2\pi \widetilde{z}}\int_{\mathbb{S}_1}\frac{dz'}{2\pi z'}\chi_\omega({z'}^{-1})\omega(\widetilde{z})p\left(z'\widetilde{z}\right)=p\left(1\right),
\end{equation}
which holds for any polynomial $p$ of order $N-1$. Indeed, for a monomial $p(z')={z'}^j$ we create a factor $1/\mathcal{S}\omega(j)$ from the $z'$-integral and a factor $\mathcal{S}\omega(j)$ in the $\widetilde{z}$-integral which obviously cancel. Thence, it is
\begin{equation}
\begin{split}
\int_{\mathbb{S}_1}\frac{d\widetilde{z}}{2\pi \widetilde{z}} P_a(\widetilde{z})Q_b(\widetilde{z})=&\int_{\mathbb{S}_1}\frac{d\widetilde{z}}{2\pi \widetilde{z}}\int_{\mathbb{S}_1}\frac{dz'}{2\pi z'}\chi_\omega({z'}^{-1})\omega(\widetilde{z})\prod_{\substack{l=1,\ldots,N\\l\neq a+1}}\frac{x_{b+1}z'\widetilde{z}-x_l}{x_{a+1}-x_l}=\delta_{ab}.
\end{split}
\end{equation}
The last equality sign is evident because the quotient $\prod_{\substack{l=1,\ldots,N\\l\neq a+1}}(x_{b+1}-x_l)/(x_{a+1}-x_l)$ vanishes whenever $l=b+1$.

For the kernel~\eqref{kernel-cyc.Pol.fixed}, we rescale the $z'$-integral in the definition of $P_j$ by $R>1$. This is essential so that when carrying out the $z'_2$-integral in~\eqref{kernel-cyc.Pol.fixed} by the residue theorem we only pick up the contributions at the $N$ poles $x_1,\ldots,x_N$. Each pole yields one summand $P_j(z_1)Q_j(z_2)$ as can be readily checked. This concludes the proof.
\hfill$\square$

\subsection{Eigenvalue Statistics of a Product comprising a Cyclic Polynomial Ensemble}\label{sec:EVstat.polynomial}

At last we consider the case from Theorem~\ref{thm:jpdfs} part (2) where $U_1\in\U(N)$ is drawn from a polynomial ensemble. We will anew make use of the polynomial~\eqref{chi-poly} when answering the question about the eigenvalue statistics at finite matrix dimension.

\begin{proposition}[Kernel of a Cyclic P\'olya Ensemble times a Cyclic Polynomial Ensemble]\label{prop:cyc.Pol.polynomial}\
	
	Let us consider the setting of Theorem~\ref{thm:jpdfs}.2 and let $\{\widetilde{P}_j,\widetilde{Q}_j\}_{j=0,\ldots,N-1}$ be a bi-orthonormal pair of functions of the cyclic polynomial random matrix $U_1\in\U(N)$. The bi-orthonormal pair of functions $\{(P_{j},Q_j)\}_{j=0,\ldots,N-1}$ for the product matrix $U=U_1U_2$ is then given by
	\begin{equation}\label{biorth-cyc.Pol.polynomial}
	\begin{split}
	P_j(z_1)=\chi_\omega\ast\widetilde{P}_j(z_1)=\int_{\mathbb{S}_1}\frac{dz'_1}{2\pi z'_1}\chi_\omega(z'_1)\widetilde{P}_j\left(\frac{z_1}{z'_1}\right),\quad Q_j(z_2)=\omega\ast\widetilde{Q}_j(z_2)=\int_{\mathbb{S}_1}\frac{dz'_2}{2\pi z'_2}\omega(z'_2)\widetilde{Q}_j\left(\frac{z_2}{z'_2}\right)
	\end{split}
	\end{equation}
	for $j=0,\ldots,N-1$ and the corresponding kernel has the following relation to the kernel $\widetilde{K}_N$ corresponding to $U_1$:
	\begin{equation}\label{kernel-cyc.Pol.polynomial}
	\begin{split}
	K_N(z_1,z_2)=&\int_{\mathbb{S}_1}\frac{dz'_1}{2\pi z'_1}\int_{\mathbb{S}_1}\frac{dz'_2}{2\pi z'_2}\chi_\omega(z'_1)\omega(z'_2)\widetilde{K}_N\left(\frac{z_1}{z'_1},\frac{z_2}{z'_2}\right).
	\end{split}
	\end{equation}
\end{proposition}

{\bf Proof of Proposition~\ref{prop:cyc.Pol.polynomial}:}

The functions $Q_j$ are inside the span of $\{\omega\ast w_j\}_{j=0,\ldots,N-1}$ because the convolution on $\mathbb{S}_1$ is linear and the functions $\widetilde{Q}_j$ are a basis of the span of $\{ w_j\}_{j=0,\ldots,N-1}$. Their linear independence can be checked by applying the Fourier transform on $Q_j=\omega\ast\widetilde{Q}_j$ and exploiting the fact that at least $N$ frequencies, namely $s=0,\ldots,N-1$, $\mathcal{S}(s)$ is invertible.

The bi-orthonormality of the pair of functions is again a direct consequence of~\eqref{poly-int-ident} as we have
\begin{equation}
\begin{split}
\int_{\mathbb{S}_1}\frac{d\widetilde{z}}{2\pi \widetilde{z}} P_a(\widetilde{z})Q_b(\widetilde{z})=&\int_{\mathbb{S}_1}\frac{d\widetilde{z}}{2\pi \widetilde{z}} \int_{\mathbb{S}_1}\frac{dz'_1}{2\pi z'_1}\int_{\mathbb{S}_1}\frac{dz'_2}{2\pi z'_2}\chi_\omega(z'_1)\omega(z'_2)\widetilde{P}_a\left(\frac{z'_2\widetilde{z}}{z'_1}\right)\widetilde{Q}_b\left(\widetilde{z}\right)\\
=&\int_{\mathbb{S}_1}\frac{d\widetilde{z}}{2\pi \widetilde{z}}\widetilde{P}_a\left(\widetilde{z}\right)\widetilde{Q}_b\left(\widetilde{z}\right)=\delta_{ab},
\end{split}
\end{equation}
where eventually have employed the bi-orthonormality of $\{\widetilde{P}_j,\widetilde{Q}_j\}_{j=0,\ldots,N-1}$.

The formula of the kernel~\eqref{kernel-cyc.Pol.polynomial} is obtained by switching the two integrals over $z'_1$ and $z'_2$ by the sum over the index $j=0,\ldots,N-1$ in~\eqref{k-point}. This is allowed as the sum is one over finite summands in the integrands are absolutely integrable because they consist of polynomials in one of the two integration arguments and of a convolution of a linear combination of $L^1$ functions in the second argument. This finishes the proof of our claims.
\hfill$\square$

\section*{Acknowledgments}

This work is part of a research program supported by the Australian Research Council (ARC) through the ARC Centre of Excellence for Mathematical and Statistical frontiers (ACEMS). Related to this, JZ acknowledges the support of a Melbourne postgraduate award and an ACEMS top up scholarship. Furthermore, PJF and MK acknowledge support from the ARC grant DP210102887. SL acknowledges the support from the National Nature Science Foundation of China (No. 12175155). We are very grateful to the referees for the detailed reports and valuable comments.

\

\end{document}